\documentclass[12pt, twoside]{article}

\usepackage{footmisc}
\usepackage{times}
\usepackage{graphics}
\usepackage{theorem}
\usepackage{graphicx}
\usepackage{amsmath}
\usepackage{latexsym}
\usepackage{amssymb}
\makeatletter
\def\eqnarray{\stepcounter{equation}\let\@currentlabel=\theequation
\global\@eqnswtrue
\global\@eqcnt\z@\tabskip\@centering\let\\=\@eqncr
$$\halign to \displaywidth\bgroup\@eqnsel\hskip\@centering
  $\displaystyle\tabskip\z@{##}$&\global\@eqcnt\@ne 
  \hfil$\;{##}\;$\hfil
  &\global\@eqcnt\tw@ $\displaystyle\tabskip\z@{##}$\hfil 
   \tabskip\@centering&\llap{##}\tabskip\z@\cr}
\makeatother
\def\thesection{\Roman{section}} % section number in Roman
\theorembodyfont{\itshape}
\newtheorem{thm}{Theorem}[section]
\newtheorem{lem}{Lemma}[section]

\newtheorem{coro}{Corollary}[section]

\theorembodyfont{\rmfamily}

\newtheorem{rem}{Remark}[section]{\bf}{\rm}
{\bf}{\rm}
\newtheorem{assumpt}{Assumption}[section]{\bf}{\it}
{\bf}{\rm}
\newcommand{\E}{\mathrm{E}}
\newcommand{\K}{\mathcal{K}}
\newcommand{\M}{\mathcal{M}}
\newcommand{\Z}{\mathcal{Z}}

\newcommand{\vc}[1]{\mbox{\boldmath $#1$}}
\newcommand{\svc}[1]{\mbox{\boldmath $\scriptstyle #1$}}

\newcommand{\wtilde}[1]{\widetilde{#1}}
\newcommand{\down}[2]{\smash{\lower#1\hbox{#2}}}
\newcommand{\up}[2]{\smash{\lower-#1\hbox{#2}}}
\newcommand{\qed}{\hspace*{\fill}$\Box$}

\makeatother

\begin{document}\thispagestyle{plain} 

\hfill
%{\small Last update date: \today}
%Submitted to STOCHASTIC MODELS, \today.

%{\it Preprint submitted to Stochastic Analysis and Applications \hfill \today}
%Submitted to STOCHASTIC MODELS, \today.

{\Large{\bf
\begin{center}
Analysis and Computation of the Joint Queue Length Distribution in a FIFO Single-Server Queue with 
Multiple Batch Markovian Arrival Streams%
\footnote[1]{This paper is published in Stochastic Models, vol. 19, no. 3, pp. 349-381, 2003.}
\end{center}
}
}

\begin{center}
{\bf 
{\bf Hiroyuki Masuyama%
\footnote[2]{Dept. of Systems Science,
Graduate School of Informatics, Kyoto University; \\\qquad\,{}%
E-mail: masuyama@sys.i.kyoto-u.ac.jp} 
and 
Tetsuya Takine%
\footnote[3]{Dept. of Information and Communications Technology, Graduate School of Engineering, Osaka University; \\\qquad\,{}%
E-mail: takine@comm.eng.osaka-u.ac.jp}
} }

\bigskip
\medskip

{\small
\textbf{Abstract}

\medskip

\begin{tabular}{p{0.85\textwidth}}
This paper considers a work-conserving FIFO single-server queue with
multiple batch Markovian arrival streams governed by a continuous-time
finite-state Markov chain.  A particular feature of this queue is that
service time distributions of customers may be different for different
arrival streams. After briefly discussing the actual waiting time
distributions of customers from respective arrival streams, we derive
a formula for the vector generating function of the
time-average joint queue length distribution in terms of the virtual
waiting time distribution.  Further assuming the discrete phase-type 
batch size distributions, we develop a numerically feasible procedure to
compute the joint queue length distribution. Some numerical examples
are provided also.
\end{tabular}
}
% \samllsize ends
\end{center}

\begin{center}
\begin{tabular}{p{0.90\textwidth}}
{\small

\medskip

}%\samllsize ends
\end{tabular}

\end{center}

\section{Introduction}

In this paper, we study the joint queue length distribution in a
stationary work-conserving FIFO single-server queue fed by multiple
batch arrival streams governed by a continuous-time finite-state
Markov chain. A particular feature of this queue is that service time
distributions of customers may be different for different arrival
streams.

Single-server queues with Markovian arrival streams have been
extensively studied for last two decades. At present, the most popular
Markovian arrival stream is MAP (Markovian arrival process) introduced
in \cite{Luca90}. MAP is a class of semi-Markovian arrival processes
including Markov modulated Poisson processes and phase-type renewal
processes as special cases. After introducing MAP, some extensions
have been made. One is batch MAP \cite{Luca91} that allows batch
arrivals and the other is marked MAP \cite{Asmu93, He96,He01} that explicitly
represents possibly correlated multiple Markovian arrival streams.
The arrival process in this paper has these two features, i.e., batch
marked MAP.

Most of previous works on FIFO single-server queues with Markovian
arrival streams assume that service times of all customers are
independent and identically distributed (i.i.d.) according to a common
distribution function. As a result, the bivariate process of the total
number of customers and the state of the Markov chain that governs the
arrival process immediately after departures forms a Markov chain of
M/G/1 type and the steady-state solution can be computed by well-known
M/G/1 paradigm \cite{Neut89}.

On the other hand, if service time distributions of customers from
respective arrival streams are different from one another, the
bivariate process does not have the Markov property \cite{Taki02},
except for queues with a superposition of independent Poisson
streams. Thus the queue length analysis of such a queue is not
straightforward. Note, however, that the virtual waiting time process
in such a queue is characterized by a bivariate Markov process
\cite{Asmu91,Taki94a,Taki96a,Zhu91}, and algorithmic solution methods
are known in the literature
\cite{Taki94a,Taki96a}.

Recently, a new approach was developed to characterize the joint queue
length distribution in FIFO queues with marked MAP having different
service time distributions \cite{Taki01b,Taki02}. In these works, the
invariant relationship of the joint queue length distributions at a
random point in time and at departures was obtained and from this, the
distributional form of Little's law was established in
\cite{Taki01b}. Further, based on the latter, an algorithmic solution
method was developed \cite{Taki01b,Taki02}. Related works are found in
\cite{Mach99,Regt86}. See \cite{Taki01a} for a survey of those
developments.

The results in this paper are considered as an extension of those in
\cite{Taki02}, allowing batch arrivals in each arrival stream. Note
here that the distributional form of Little's law does not hold for
FIFO queues with batch arrivals. Therefore our starting point in
analyzing the time-average joint queue length distribution is the
invariant relationship of the joint queue length distributions at a
random point in time and at departures in \cite{Taki01b}. By doing so,
the problem is reduced to find the joint queue length distributions at
departures of customers from respective arrival streams.

As you will see, the joint queue length distribution at departures in
the FIFO queue is closely related to the virtual waiting time
distribution that is readily obtained with the known results. Using
these facts, we derive a general formula for the stationary joint
queue length distribution at departures in terms of the sojourn time
distribution. Further, assuming discrete phase-type batch size
distributions, we derive recursions to compute the joint queue length
distribution.

The above outline is similar to the single arrival case in
\cite{Taki02}. However, the implementation of some of those recursions
is not trivial, because we have to determine several truncation and
stopping criteria, which are due to batch arrivals, and their
straightforward implementation would require very huge memory space
and time-consuming. In this paper, assuming discrete phase-type batch
size distributions, we propose a numerically feasible procedure to
compute those recursions, while ensuring the numerical accuracy in the
final result. This is the main contribution of this paper.  Note that
our procedure is applicable to the FIFO BMAP/G/1 queue with i.i.d.\
services. too, when the batch size distribution follows a discrete
phase-type distribution.

The rest of this paper is divided into six sections. In section
\ref{sec:2}, the mathematical model is described. In section
\ref{sec:3}, we briefly discuss the virtual and actual waiting time
distributions. In section \ref{sec:4}, we first derive a general
formula for the joint queue length distribution, and assuming the
discrete phase-type batch sizes, we show recursive formulas to compute
the joint queue length distribution. In section \ref{sec:5}, the
implementation of the recursions is discussed. In section \ref{sec:6},
we discuss the efficiency of our algorithm and the qualitative
behavior of the queue length through some numerical examples. Finally,
concluding remarks are provided in section \ref{sec:7}.  Throughout
the paper, we denote matrices and vectors by bold capital letters and
bold small letters, respectively.

\section{Model}\label{sec:2}

We consider a work-conserving FIFO single-server queue fed by $K$
arrival streams. We call customers arriving from the $k$th
($k=1,\ldots,K$) arrival stream class $k$ customers. Let $\K$ denote a
set of class indices, i.e., $\K = \{1,2,\dots,K \}$.

Customer arrivals are governed by a continuous-time Markov chain,
which is called the underlying Markov chain hereafter. The underlying
Markov chain has a finite state space ${\cal M}=\{1,\ldots,M\}$ and it
is assumed to be irreducible. The underlying Markov chain stays in
state $i \in \M$ for an exponential interval of time with mean
$\mu_i^{-1}$. When the sojourn time in state $i$ has elapsed, with
probability $\sigma_{i,j}(0)$ ($j \in \M$, $j \neq i$), the underlying
Markov chain changes its state to state $j$ without arrivals. Also,
with probability $\sigma_{k,i,j}(n)$ ($k \in \K$, $n = 1,2,\ldots$),
the underlying Markov chain changes its state to state $j$ and $n$
customers of class $k$ arrive simultaneously. For convenience, let
$\sigma_{i,i}(0) = 0$ for all $i \in \M$. Then
\[
\sum_{j \in \M} \left(\sigma_{i,j}(0)
+  \sum_{k \in \K} \sum_{n=1}^{\infty} \sigma_{k,i,j}(n)\right) = 1,
\]%
\noindent
for all $i \in \M$. We assume that service times of class $k$ ($k \in
\K$) customers are i.i.d.\ according to a distribution function
$H_k(x)$ with finite mean $h_k$.

We now introduce some notations to describe the above arrival process. 
Let $\vc{C}$ denote an $M \times M$ matrix whose $(i,j)$th $(i,j \in
\M)$ element $C_{i,j}$ is given by
\[
C_{i,j} = 
\left\{
\begin{array}{ll}
-\mu_i, 
& \mbox{if\ } i = j,  \\
\sigma_{i,j}(0) \mu_i, 
& \mbox{otherwise}.
\end{array}
\right.
\]%
\noindent
Further, for $k \in \K$, we define $\vc{D}_k(n)$ ($n = 1,2,\ldots$) as
an $M \times M$ matrix whose $(i,j)$th $(i,j \in \M)$ element
$D_{k,i,j}(n)$ is given by
\[
D_{k,i,j}(n) = \sigma_{k,i,j}(n) \mu_i.
\]%
\noindent
Thus the counting process of arrivals is characterized by the set of
matrices $(\vc{C},\vc {D}_1(n_1),\dots, \vc{D}_K(n_K))$. Roughly
speaking, customers arrive in the following way. When a state
transition driven by $\vc{D}_k(n)$ occurs, $n$ customers of class $k$
arrive simultaneously. On the other hand, when a state transition
driven by $\vc{C}$ occurs, no customers arrive.

We define $\vc{D}_k$ $(k \in \K)$ and $\vc{D}$ as
\[
\vc{D}_k = \sum_{n=1}^{\infty} \vc{D}_k(n),
\qquad 
\vc{D} = \sum_{k \in \K} \vc{D}_k, 
\]%
\noindent
respectively. Note that the infinitesimal generator of the underlying
Markov chain is given by $\vc{C} + \vc{D}$. Note also that $(\vc{C} +
\vc{D}) \vc{e} = \vc{0}$, where $\vc{e}$ denotes a column 
vector whose elements are all equal to one. We denote, by $\vc{\pi}$,
the stationary probability vector of the underlying Markov chain and
therefore $\vc{\pi}$ satisfies $\vc{\pi}(\vc{C} + \vc{D}) = \vc{0}$
and $\vc{\pi}\vc{e} = 1$. Because of the finite state space $\M$
and the irreducibility of the underlying Markov chain, $\vc{\pi}$ is
uniquely determined.

We define $\lambda_k$ $(k \in \K)$ as 
\[
\lambda_k = \sum_{n=1}^{\infty} n \vc{\pi} \vc{D}_k(n)\vc{e}.
\]%
\noindent
Note that $\lambda_k$ denotes the arrival rate of class $k$ customers,
i.e., the mean number of class $k$ customers arriving in a unit time
in steady state. We assume that at least one element of $\vc{D}_k$ ($k
\in \K$) is positive, so that $\lambda_k > 0$ for all $k \in \K$. Let
$\rho_k$ denote the utilization factor of class $k$ customers,\ i.e.,
\[
\rho_k = \lambda_k h_ k, \qquad k \in \K.
\]%
\noindent
Furthermore, we denote the overall arrival rate by $\lambda = \sum_{k
\in \K} \lambda_k$ and the overall utilization factor by $\rho =
\sum_{k \in \K} \rho_k$. In the remainder of this paper, we assume
that $\rho < 1$, which ensures that all customers arriving to the
system are eventually served \cite{Loyn62}.

\section{Waiting Time Distribution}\label{sec:3}

In this section, we consider the stationary distribution of the actual
waiting time. To do so, we first consider the virtual waiting time
that is equivalent to the amount of work in system. Let $V$ denote a
generic random variable representing the stationary amount of work in
system (i.e., the total amount of unfinished services of all customers
in the system). Also let $S$ denote a generic random variable
representing the state of the underlying Markov chain in steady state. 
We then define $\vc{v}(x)$ as a $1 \times M$ vector whose $j$th
element represents $\Pr\left[V \le x, S = j \right]$.  The
Laplace-Stielties transforms (LSTs) of $H_k(x)$ and $\vc{v}(x)$ are
denoted by $H_k^{\ast}(s)$ and $\vc{v}^{\ast}(s)$, respectively.  

We define $\vc{D}(x)$ as
\[
\vc{D}(x) = \sum_{k \in \K} \sum_{n=1}^{\infty} 
\vc{D}_k(n) H_k^{(n)}(x),
\qquad x \ge 0,
\]%
\noindent
where $H_k^{(1)}(x) = H_k(x)$ and $H_k^{(n)}(x)$ $(n=2,3,\dots)$
denotes the $n$-fold convolution of $H_k(x)$ with itself. 
Let $\vc{Q}$ denote an $M \times M$ matrix that represents the
infinitesimal generator of the underlying Markov chain obtained by
excising the busy periods \cite{Taki94a}. Note that $\vc{Q}$ satisfies
\[
\vc{Q} = \vc{C} + \int_0^{\infty}d\vc{D}(x) \exp(\vc{Q} x).
\]%
\noindent
Let $\vc{\kappa}$ denote a $1 \times M$ vector that satisfies
\[
\vc{\kappa}\vc{Q} = \vc{0}, \qquad \vc{\kappa}\vc{e} = 1.
\]%
\noindent
Applying the results in \cite{Taki94a} to our model, 
we obtain the following theorem.
\begin{thm}[\cite{Taki94a}]
$\vc{v}(0)$ is given by
\[
\vc{v}(0) = (1 - \rho)\vc{\kappa}.
\]%
\noindent
Furthermore, the LST $\vc{v}^{\ast}(s)$ of $\vc{v}(x)$ satisfies
\begin{equation}
\vc{v}^{\ast}(s)
\left[s \vc{I} + \vc{C} 
+ \vc{D}^{\ast}(s) \right] 
= s(1 - \rho)\vc{\kappa}, 
\qquad \mbox{\rm Re}(s) > 0,
\label{eqn:a-11}
\end{equation}%
\noindent
where $\vc{D}^{\ast}(s)$ denotes the LST of $\vc{D}(x)$:
\begin{eqnarray}
\vc{D}^{\ast}(s)
= \int_0^{\infty} e^{-sx}d\vc{D}(x)
= \sum_{k \in \K} \sum_{n=1}^{\infty} \vc{D}_k(n) 
\{ H_k^{\ast}(s) \}^{n}.
\label{eqn:i-14}
\end{eqnarray}%
\end{thm}

We now consider the actual waiting time of class $k$ customers in
steady state. We define $W_k(n;m)$ as a generic random variable
representing the actual waiting time of a randomly chosen class $k$
customer who is a member of a batch of size $n$ and the $m$th served
customer among members of the same batch. Let
$S^{(\mathrm{A}_{k})}(n)$ denote a generic random variable
representing the state of the underlying Markov chain immediately
after class $k$ batches of size $n$ arrive. With those, we define
$\vc{w}_k(x\, |\, n;m)$ as a $1 \times M$ vector whose $j$th element
represents $\Pr [W_k(n;m) \le x, S^{(\mathrm{A}_{k})}(n) = j]$. Note
that
$\vc{w}_k(x \, | \, n;1)$ $(n \ge 1)$ is given by \cite{Taki94a}
\begin{equation}
\vc{w}_k(x\, | \, n;1)
= {\vc{v}(x)\vc{D}_k(n) \over \vc{\pi}\vc{D}_k(n) \vc{e}},
\label{eqn:a-02}
\end{equation}%
\noindent
and for $n \ge 2$ and $m=2,3,\dots,n$,
\begin{equation}
\vc{w}_k(x\, |\, n;m) 
= \int_0^x\vc{w}_k(x-y\, |\,n;1) dH_k^{(m-1)}(y),
\qquad x \ge 0.
\label{eqn:f-07}
\end{equation}%

Let $W_k$ and $S^{(\mathrm{A}_{k})}$ denote generic random variables
representing the actual waiting time of class $k$ customers and the
state of the underlying Markov chain immediately after arrivals of
class $k$ batches, respectively. We then define $\vc{w}_k(x)$ as a $1
\times M$ vector whose $j$th element represents $ \Pr [W_k \le x,
S^{(\mathrm{A}_{k})} = j]$. Because a randomly chosen customer of
class $k$ is a member of a batch of size $n$ with probability
$n\vc{\pi}\vc{D}_k(n) \vc{e}/\lambda_k$, we have
\begin{equation}
\vc{w}_k(x) 
= \sum_{n=1}^{\infty} {n \vc{\pi}\vc{D}_k(n)\vc{e} \over \lambda_k}
\cdot {1 \over n} \sum_{m=1}^n \vc{w}_k(x\, | \, n;m).
\label{eqn:a-14}
\end{equation}%
\noindent
Let $\vc{w}_k^{\ast}(s)$ denote the LST of $\vc{w}_k(x)$.  From
(\ref{eqn:a-02})--(\ref{eqn:a-14}), we have
\[
\vc{w}_k^{\ast}(s) 
= \sum_{n=1}^{\infty} 
{\vc{v}^{\ast}(s)\vc{D}_k(n) \over \lambda_k}
\sum_{m=1}^n \{H_k^{\ast}(s)\}^{m-1},
\qquad \mbox{Re}(s) > 0.
\]%
\noindent
Thus we obtain the following theorem.
\begin{thm}\label{thm:b-01}
$\vc{w}_k^{\ast}(s)$ ($k \in \K$) is given by
\[
\vc{w}_k^{\ast}(s) = {\vc{v}^{\ast}(s)\left(\vc{D}_k -
\vc{D}_k^{\ast}(H_k^{\ast}(s))\right) \over \lambda_k \left( 1 -
H_k^{\ast}(s) \right)}, \qquad \mbox{\rm Re}(s) > 0,
\]%
\noindent
where
\begin{equation}
\vc{D}_k^{\ast}(z_k) = \sum_{n=1}^{\infty} z_k^{n} \vc{D}_k(n).
\label{eqn:b-14}
\end{equation}%
\end{thm}

\section{Joint Queue Length Distribution}\label{sec:4}

This section considers the joint queue length distribution. In
subsection \ref{subsec:A}, we apply a general relationship between the
time-average queue length distribution and the queue length
distributions at departures of customers of respective classes
\cite{Taki01b} to our specific queue. Then the problem is reduced to
characterize the joint queue length distributions at departures of
respective classes, which is discussed in subsection
\ref{subsec:B}. Finally in subsection \ref{subsec:C}, assuming
discrete phase-type batch size distributions, we derive recursions for
some quantities required in computing the joint queue length
distribution.

\subsection{Relationship in the joint queue length distributions}\label{subsec:A}

Let $N_k$ $(k \in \K)$ denote a generic random variable representing
the number of class $k$ customers in steady state. We define
$\vc{p}(n_1, \dots, n_K)$ as a $1 \times M$ vector whose $j$th element
represents $\Pr[N_1=n_1, \ldots$, $N_K=n_K, S = j]$. For simplicity,
let $\vc{n}$ and $\vc{z}$ denote a $1 \times K$ nonnegative integer
vector $(n_1,\dots,n_K)$ and a $1 \times K$ complex vector
$(z_1,\dots,z_K)$, respectively. Further we define $\Z$ as
\[
\Z = \{ (n_1,\dots,n_K) \,
;\, n_k = 0,1,\dots, \mbox{ for all } k \in \K \}.
\]%
\noindent
We then define $\vc{p}^{\ast}(\vc{z})$ as
\[
\vc{p}^{\ast}(\vc{z})
= \sum_{\svc{n} \in \Z} z_1^{n_1} \cdots
z_K^{n_K} \vc{p}(\vc{n}), 
\qquad |z_k| \leq 1  \mbox{ for all } k \in \K.
\]%
\noindent
Note that $\vc{p}^{\ast}(\vc{z})$ denotes the vector generating
function of the joint queue length distribution in steady state.

Let $N_{\nu}^{(\mathrm{D}_k)}$ and $S^{(\mathrm{D}_k)}$ ($k,\nu \in
\K$) denote generic random variables representing the number of class
$\nu$ customers and the state of the underlying Markov chain,
respectively, immediately after departures of class $k$ customers in
steady state. We then define $\vc{q}_k(\vc{n})$ ($k \in \K$, $\vc{n}
\in \Z$) as a $1 \times M$ vector whose $j$th element represents
$\Pr[N_1^{(\mathrm{D}_k)} = n_1, \ldots, N_K^{(\mathrm{D}_k)} = n_K,
S^{(\mathrm{D}_k)} = j]$. Further we define $\vc{q}_k^{\ast}(\vc{z})$
$(k \in \K)$ as
\[
\vc{q}_k^{\ast}(\vc{z}) 
= \sum_{\svc{n} \in \Z} z_1^{n_1} \cdots
z_K^{n_K} \vc{q}_k(\vc{n}),
\qquad |z_k| \leq 1  \mbox{ for all } k \in \K.
\]%
\noindent
Note that $\vc{q}_k^{\ast}(\vc{z})$ denotes the vector generating
function of the joint queue length distribution immediately after
departures of class $k$ customers. Thus, applying Theorem 1 in
\cite{Taki01b} to our model, we have the following theorem.
\begin{thm}[\cite{Taki01b}]\label{thm:a-02}
$\vc{p}^{\ast}(\vc{z})$ and $\vc{q}_k^{\ast}(\vc{z})$ are related by
\begin{equation}
\vc{p}^{\ast}(\vc{z}) \left[ \vc{C} 
+ \sum_{k \in \K} \vc{D}_k^{\ast}(z_k) \right]
= \sum_{k \in \K} \lambda_k(z_k - 1)\vc{q}_k^{\ast}(\vc{z}),
\label{ean:b-05}
\end{equation}%
\noindent
where $\vc{D}_k^{\ast}(z_k)$ is given in (\ref{eqn:b-14}).
\end{thm}

Further, comparing the coefficient vectors of $z_1^{n_1} \cdots
z_K^{n_K}$ on both sides of (\ref{ean:b-05}), we obtain the following
result.
\begin{coro}\label{coro:a-04}
The $\vc{p}(\vc{n})$ $(\vc{n} \in \Z)$ is recursively determined by
\begin{eqnarray*}
\vc{p}(\vc{0}) 
&=& \sum_{k \in \K} \lambda_k \vc{q}_k(\vc{0})(-\vc{C})^{-1},
\\
\vc{p}(\vc{n})
&=& \sum_{k \in \K} 
\bigg[ \lambda_k 
\left( \vc{q}_k(\vc{n}) - \vc{q}_k(\vc{n}-\vc{e}_k)
\right)
 + \sum_{m_k=1}^{n_k} 
\vc{p}(\vc{n}-m_k\vc{e}_k) \vc{D}_k(m_k) 
\bigg] (-\vc{C})^{-1}, 
\qquad \vc{n} \in \Z^+ ,
\end{eqnarray*}%
\noindent
where $\Z^{+} = \Z - \{\vc{0}\}$, $\vc{q}_k(\vc{n}) = \vc{0}$ 
for $\vc{n}\ \ooalign{\hfill $\in$\hfill\crcr\hfill /\hfill}\ \Z$
and $\vc{e}_k$ $(k \in \K)$ denotes the $k$th unit vector:
\[
\vc{e}_k = (0,\ldots,0,\begin{array}[t]{@{}c@{}}1\\ 
\hbox to 0mm{\hss $k$\mbox{\rm th}\hss}
\end{array},0,\ldots,0).
\]%
\end{coro}

\begin{rem}
Throughout the paper, the empty sum is defined as zero.
\end{rem}

\subsection{Joint queue length distribution immediately after departures}\label{subsec:B}

In this subsection, we consider the vector generating function of the
joint queue length distribution immediately after departures of each
class. We denote, by $\mathrm{C}_k(n;m)$ ($k \in \K$, $n=1,2,\ldots$,
$m=1,2,\dots,n$), a randomly chosen class $k$ customer who is a member
of a batch of size $n$ and the $m$th served customer among members of
the same batch.  Let $N_{\nu}^{(\mathrm{D}_k)}(n;m)$ and
$S^{(\mathrm{D}_k)}(n;m)$ ($k,\nu \in \K$, $n=1,2,\ldots$,
$m=1,2,\dots,n$) denote generic random variables representing the
number of class $\nu$ customers and the state of the underlying Markov
chain, respectively, immediately after the departure of customer
$\mathrm{C}_k(n;m)$ in steady state. We then define
$\vc{q}_k^{\ast}(\vc{z}\, | \, n;m)$ ($k \in \K$, $n=1,2,\ldots$,
$m=1,2,\dots,n$) as a $1 \times M$ vector whose $j$th element
represents
\[
\E \left[\prod_{\nu \in \K}
z_{\nu}^{N_{\nu}^{(\mathrm{D}_k)}(n;m)} 
1\{S^{(\mathrm{D}_k)}(n;m) = j\} \right],
\]%
\noindent
where $1\{\chi\}$ denotes an indicator function of event $\chi$.
Because a randomly chosen customer of class $k$ is a member of a batch
of size $n$ with probability $n\vc{\pi}\vc{D}_k(n) \vc{e}/\lambda_k$,
we have
\begin{equation}
\vc{q}_k^{\ast}(\vc{z}) 
= \sum_{n=1}^{\infty}
{n\vc{\pi}\vc{D}_k(n) \vc{e} \over \lambda_k}
{1 \over n} \sum_{m=1}^{n} \vc{q}_k^{\ast}(\vc{z}\, | \, n;m).
\label{eqn:a-24}
\end{equation}%
\noindent
In what follows, we consider $\vc{q}_k^{\ast}(\vc{z}\, | \, n;m)$.

We define $\overline{W}_k(n;m)$ ($k \in \K$, $n=1,2,\ldots$,
$m=1,2,\dots,n$) as a generic random variable representing the sojourn
time of customer $\mathrm{C}_k(n;m)$. Note here that
\[
\overline{W}_k(n;m) = W_k(n; 1)+ H_{k,1} + \cdots + H_{k,m},
\]%
\noindent
where $W_k(n; 1)$ denotes the actual waiting time of customer
$\mathrm{C}_k(n;1)$, and $H_{k,l}$ ($l=1,2,\dots,m$) denotes the
service time of customer $\mathrm{C}_k(n;l)$. By definition, $W_k(n;
1)$ depends only on the past history up to the arrival instant of a
batch including customer $\mathrm{C}_k(n;m)$. On the other hand, the
number of customers in the system immediately after the departure of
customer $\mathrm{C}_k(n;m)$ is equal to the sum of the $n-m$
customers in the same batch and customers who arrived during the
sojourn time of customer $\mathrm{C}_k(n;m)$. Note here that the
latter is conditionally independent of the past history given the
length of the sojourn time and the state of the underlying Markov
chain immediately after the arrival of the batch. Thus we have
\begin{eqnarray}
\vc{q}_k^{\ast}(\vc{z}\, | \, n;m)
&=& z_k^{n-m}
\int_0^{\infty} d\vc{w}_k(x \, | \, n;1)
\vc{N}^{\ast}(x, \vc{z}) 
 \left[ \int_0^{\infty} dH_k(y)  
\vc{N}^{\ast}(y, \vc{z})\right]^m,
\label{eqn:a-01}
\end{eqnarray}%
\noindent
where
\begin{equation}
\vc{N}^{\ast}(x, \vc{z})
= 
\exp\left[ \left(\vc{C} + \sum_{k \in \K} \vc{D}_k^{\ast}(z_k) 
\right) x\right].
\label{eqn:i-17}
\end{equation}%
\begin{thm}
The vector generating function $\vc{q}_k^{\ast}(\vc{z})$ ($k \in \K$)
of the joint queue length distribution immediately after departures of
class $k$ customers in the steady state is given by
\begin{eqnarray}
\vc{q}_k^{\ast}(\vc{z})
&=& {1 \over \lambda_k}
\sum_{m=1}^{\infty} \sum_{l=0}^{\infty} z_k^l
\int_0^{\infty} d\vc{v}(x)\vc{D}_k(m+l)
\vc{N}^{\ast}(x, \vc{z})
\left[ \int_0^{\infty} dH_k(y)  
\vc{N}^{\ast}(y, \vc{z}) \right]^m.
\label{eqn:a-07}
\end{eqnarray}%
\end{thm}

\noindent
\textbf{Proof.}
Using (\ref{eqn:a-02}), (\ref{eqn:a-24}) and (\ref{eqn:a-01}), we have
\begin{eqnarray*}
\vc{q}_k^{\ast}(\vc{z})
&=& \sum_{n=1}^{\infty}
{\vc{\pi} \vc{D}_k(n) \vc{e} \over \lambda_k}
\sum_{m=1}^n z_k^{n-m} \int_0^{\infty}
{d\vc{v}(x)\vc{D}_k(n) \over \vc{\pi} \vc{D}_k(n) \vc{e}}
\vc{N}^{\ast}(x, \vc{z})
\left[ \int_0^{\infty} dH_k(y)  
\vc{N}^{\ast}(y, \vc{z}) \right]^m
\\
&=& {1\over \lambda_k}\sum_{n=1}^{\infty}
\sum_{m=1}^n z_k^{n-m} \int_0^{\infty}
d\vc{v}(x)\vc{D}_k(n)
\vc{N}^{\ast}(x, \vc{z})
\left[ \int_0^{\infty} dH_k(y)  
\vc{N}^{\ast}(y, \vc{z}) \right]^m,
\end{eqnarray*}%
\noindent
from which (\ref{eqn:a-07}) follows. \qed

\subsection{Recursions for discrete phase-type batch sizes}\label{subsec:C}

In this subsection, we develop a recursive formula to compute the
vector mass function $\vc{q}_k(\vc{n})$ of the joint queue length
immediately after departures of each class under the following
assumption.
\begin{assumpt}\label{ass:a-01}
The batch size distribution of class $k$ is independent of the state
of the underlying Markov chain and follows a discrete phase-type
distribution with representation $(\vc{\alpha}_k, \vc{P}_k)$, i.e.,
\begin{eqnarray}
\vc{D}_k(n) 
&=& g_k(n)\vc{D}_k, 
\label{eqn:f-14}
\\
g_k(n) 
&=& \vc{\alpha}_k \vc{P}_k^{n-1}(\vc{I} - \vc{P}_k)\vc{e},
\qquad n=1,2,\dots,
\label{eqn:a-12}
\end{eqnarray}%
\noindent
where $\vc{\alpha}_k$ denotes a $1 \times M_k$ probability vector and
$\vc{P}_k$ denotes an $M_k \times M_k$ substochastic matrix.
\end{assumpt}

Let $\vc{I}(m)$ denote an $m \times m$ identity matrix.  When the size
of an identity matrix is clear from the context, we suppress $(m)$.
\begin{lem}
Under Assumption \ref{ass:a-01}, $\vc{q}_k^{\ast}(\vc{z})$ ($k \in
\K$) is given by
\begin{eqnarray}
\vc{q}_k^{\ast}(\vc{z}) 
&=& {1 \over \lambda_k}
\int_0^{\infty} d\vc{v}(x) \vc{D}_k \vc{N}^{\ast}(x, \vc{z}) 
\cdot 
\left( \vc{\alpha}_k \otimes \int_0^{\infty} dH_k(y)
\vc{N}^{\ast}(y, \vc{z}) \right)
\nonumber
\\
& & \qquad {} \cdot 
\left[ \vc{I} - \vc{P}_k \otimes \int_0^{\infty} dH_k(y)  
\vc{N}^{\ast}(y, \vc{z}) \right]^{-1}
\nonumber
\\
&& \qquad \quad {} \cdot
\left[
\left\{ \left(\vc{I} - z_k\vc{P}_k \right)^{-1}
(\vc{I} - \vc{P}_k) \vc{e} \right\} \otimes \vc{I}(M) \right]. 
\label{eqn:b-10}
\end{eqnarray}%
\end{lem}

\noindent
\textbf{Proof.}
Substituting (\ref{eqn:f-14}) and (\ref{eqn:a-12}) into
(\ref{eqn:a-07}) and using properties of Kronecker product:
\begin{eqnarray*}
& & a \vc{A} = a \otimes \vc{A} \mbox{ for any scalar $a$, and},
\\
& & (\vc{A}_1 \cdots \vc{A}_n) 
\otimes (\vc{B}_1 \cdots \vc{B}_n) 
\nonumber
\\
&& \qquad \quad 
= (\vc{A}_1 \otimes \vc{B}_1) 
 \cdots  (\vc{A}_n \otimes \vc{B}_n)
\mbox{ for any $n =1,2,\ldots$},
\end{eqnarray*}%
\noindent
we obtain 
\begin{eqnarray*}
\vc{q}_k^{\ast}(\vc{z}) 
&=& {1 \over \lambda_k}
\int_0^{\infty} d\vc{v}(x) \vc{D}_k 
\vc{N}^{\ast}(x, \vc{z}) 
\left\{
\vc{\alpha}_k \sum_{m=1}^{\infty} \vc{P}_k^{m-1}
\sum_{l=0}^{\infty} z_k^l \vc{P}_k^l
(\vc{I} - \vc{P}_k) \vc{e} \right\}
\nonumber
\\
&& \qquad \otimes 
\left\{
\left( \int_0^{\infty} dH_k(y) \vc{N}^{\ast}(y, \vc{z}) \right)^m
\right\}
\\
&=& {1 \over \lambda_k}
\int_0^{\infty} d\vc{v}(x) \vc{D}_k 
\vc{N}^{\ast}(x, \vc{z})  
\left\{
\vc{\alpha}_k \sum_{m=1}^{\infty} \vc{P}_k^{m-1}
\left(\vc{I} - z_k \vc{P}_k \right)^{-1}
(\vc{I} - \vc{P}_k) \vc{e} \right\}
\nonumber
\\
&& \qquad \otimes 
\left\{
\left( \int_0^{\infty} dH_k(y) \vc{N}^{\ast}(y, \vc{z}) \right)^m
\right\}
\\
&=& {1 \over \lambda_k}
\int_0^{\infty} d\vc{v}(x) \vc{D}_k \vc{N}^{\ast}(x, \vc{z})
\left( 
\vc{\alpha}_k \otimes \int_0^{\infty} dH_k(y)
\vc{N}^{\ast}(y, \vc{z}) \right) \nonumber
\\
& & \qquad {} \cdot \sum_{m=1}^{\infty} 
\left( \vc{P}_k \otimes \int_0^{\infty} dH_k(y)  
\vc{N}^{\ast}(y, \vc{z}) \right)^{m-1} 
\nonumber
\\
& & \qquad {}
\cdot 
\left[
\left\{ \left(\vc{I} - z_k\vc{P}_k \right)^{-1}
(\vc{I} - \vc{P}_k) \vc{e} \right\} \otimes \vc{I}(M) 
\right]
\\
&=& {1 \over \lambda_k}
\int_0^{\infty} d\vc{v}(x) \vc{D}_k \vc{N}^{\ast}(x, \vc{z})
\left( 
\vc{\alpha}_k \otimes \int_0^{\infty} dH_k(y)
\vc{N}^{\ast}(y, \vc{z}) \right) \nonumber
\\
& & \qquad {} \cdot 
\left[ \vc{I} - \vc{P}_k \otimes \int_0^{\infty} dH_k(y)  
\vc{N}^{\ast}(y, \vc{z}) \right]^{-1}
\nonumber
\\
& & \qquad {}
\cdot 
\left[
\left\{ \left(\vc{I} - z_k\vc{P}_k \right)^{-1}
(\vc{I} - \vc{P}_k) \vc{e} \right\} \otimes \vc{I}(M) 
\right],
\end{eqnarray*}%
\noindent
which completes the proof.
\qed

We define $\vc{v}_k(\vc{n})$ ($k \in \K$, $\vc{n} \in \Z$) as a $1
\times M$ vector satisfying
\begin{equation}
\sum_{\svc{n} \in \Z} z_1^{n_1} \cdots z_K^{n_K} 
\vc{v}_k(\vc{n}) 
= \int_0^{\infty} d\vc{v}(x) \vc{D}_k 
\vc{N}^{\ast}(x, \vc{z}).
\label{eqn:b-11}
\end{equation}%
\noindent
We also define $\vc{A}_k(\vc{n})$ and $\vc{\Gamma}_k(\vc{n})$ ($k \in
\K$, $\vc{n} \in \Z$) as $M \times M$ and $M M_k \times M M_k$
matrices satisfying
\begin{eqnarray}
\sum_{\svc{n} \in \Z} z_1^{n_1} \cdots z_K^{n_K} 
\vc{A}_k(\vc{n}) 
&=& \int_0^{\infty} dH_k(y)
\vc{N}^{\ast}(y, \vc{z}), 
\label{eqn:b-12}
\\
\sum_{\svc{n} \in \Z} z_1^{n_1} 
\cdots z_K^{n_K} \vc{\Gamma}_k(\vc{n}) 
&=& \left[ \vc{I} - \vc{P}_k \otimes \int_0^{\infty} dH_k(y)  
\vc{N}^{\ast}(y, \vc{z})\right]^{-1},
\nonumber
\\
\label{eqn:b-13}
\end{eqnarray}%
\noindent
respectively. Note here that
\begin{eqnarray*}
\left\{ \left(\vc{I} - z_k\vc{P}_k \right)^{-1}
(\vc{I} - \vc{P}_k) \vc{e} \right\} \otimes \vc{I}(M)
&=& \sum_{m=0}^{\infty} z_k^m \left\{
\vc{P}_k^m ( \vc{I} - \vc{P}_k ) \vc{e} \right\} \otimes \vc{I}(M).
\end{eqnarray*}%
\noindent
Thus (\ref{eqn:b-10}) is rewritten to be
\begin{eqnarray}
\vc{q}_k^{\ast}(\vc{z})
&=& {1 \over \lambda_k}
\sum_{\svc{n} \in \Z} z_1^{n_1} \cdots z_K^{n_K}
\sum_{m=0}^{n_k}
\sum_{\scriptstyle \svc{n}_1+\svc{n}_2 +\svc{n}_3
\atop \scriptstyle = \, \svc{n}-m\svc{e}_k }
\vc{v}_k(\vc{n}_1) 
\nonumber
\\
& & {} \cdot 
[\vc{\alpha}_k \otimes \vc{A}_k(\vc{n}_2) ] 
\vc{\Gamma}_k(\vc{n}_3)
\left[ \left\{ \vc{P}_k^m (\vc{I} - \vc{P}_k) \vc{e} \right\}
\otimes \vc{I}(M) \right], 
\nonumber
\\
\label{eqn:f-12}
\end{eqnarray}%
\noindent
where $\vc{n}_j \in \Z$ for $j=1,2,3$. Comparing coefficient vectors
of $z_1^{n_1}$ $\cdots$ $z_K^{n_K}$ on both sizes of (\ref{eqn:f-12}),
we obtain the following result.
\begin{thm}\label{thm:a-01}
Under Assumption \ref{ass:a-01}, $\vc{q}_k(\vc{n})$ ($k \in \K$,
$\vc{n} \in \Z$) is given by
\begin{eqnarray*}
\vc{q}_k(\vc{n}) 
&=& {1 \over \lambda_k}
\sum_{m=0}^{n_k}
\sum_{\scriptstyle \svc{n}_1+\svc{n}_2 +\svc{n}_3
\atop \scriptstyle = \, \svc{n}-m\svc{e}_k }
\vc{v}_k(\vc{n}_1) [ \vc{\alpha}_k \otimes \vc{A}_k(\vc{n}_2)] 
\vc{\Gamma}_k(\vc{n}_3)
\left[ \vc{P}_k^m (\vc{I} - \vc{P}_k) \vc{e} 
\otimes \vc{I}(M) \right],
\end{eqnarray*}%
\noindent
where $\vc{n}_j \in \Z$ for $j=1,2,3$.
\end{thm}
Theorem \ref{thm:a-01} implies that the computation of
$\vc{q}_k(\vc{n})$ is reduced to those of $\vc{v}_k(\vc{n})$,
$\vc{A}_k(\vc{n})$ and $\vc{\Gamma}_k(\vc{n})$, which are discussed in
the rest of this subsection.

We first consider the $\vc{A}_k(\vc{n})$. Let $\theta$ denote the
maximum absolute value of diagonal elements of $\vc{C}$.  We define
$\vc{F}_m(\vc{n})$ $(m =0,1,\ldots,\ \vc{n} \in \Z)$ as an $M
\times M$ matrix that satisfies
\begin{equation}
\sum_{\svc{n} \in \Z} z_1^{n_1} 
\cdots z_K^{n_K} \vc{F}_m(\vc{n})
=
\left[\vc{I} + \theta^{-1}
\left(\vc{C} + \sum_{k \in \K} \vc{D}_k^{\ast}(z_k)\right)
\right]^m.
\label{eqn:a-04}
\end{equation}%
\begin{lem}\label{lemma:a-01}
$\vc{A}_k(\vc{n})$ is given by
\begin{equation}
\vc{A}_k(\vc{n}) = \sum_{m=0}^{\infty} 
\gamma_k^{(m)}(\theta) \vc{F}_m(\vc{n}),
\qquad k \in \K,\ \vc{n} \in \Z ,
\label{eqn:b-22}
\end{equation}%
\noindent
where 
\begin{equation}
\gamma_k^{(m)}(\theta) = \int_0^{\infty} 
e^{-\theta y}{(\theta y)^m \over m!}
dH_k(y),
\;\, k \in \K, \ m=0,1,\ldots,
\label{eqn:a-13}
\end{equation}%
\noindent
and $\vc{F}_m(\vc{n})$'s are recursively determined by
\begin{equation}
\vc{F}_0(\vc{n}) =
\left\{
\begin{array}{ll}
\vc{I}, & \mbox{\rm if } \vc{n} = \vc{0}, \\
\vc{O}, & \mbox{\rm otherwise} ,
\end{array}
\right.
\label{eqn:add-F_0}
\end{equation}%
\noindent
and for $m=0,1,\ldots$, 
\begin{eqnarray}
\vc{F}_{m+1}(\vc{n})
&=& \vc{F}_m(\vc{n})(\vc{I} + \theta^{-1}\vc{C}) 
+ \theta^{-1}\sum_{k \in \K} \sum_{l_k=1}^{n_k}
\vc{F}_m(\vc{n} - l_k\vc{e}_k) \vc{D}_k(l_k),
\quad \vc{n} \in \Z .
\label{eqn:g-16}
\end{eqnarray}%
\end{lem}

\noindent
\textbf{Proof.} From (\ref{eqn:i-17}), (\ref{eqn:b-12}) and
(\ref{eqn:a-13}), we obtain
\begin{eqnarray}
\sum_{\svc{n} \in \Z} z_1^{n_1} \cdots z_K^{n_K} 
\vc{A}_k(\vc{n}) 
&=& 
\int_0^{\infty} dH_k(y)
\exp\left[ \left(\vc{C} + \sum_{k \in \K} \vc{D}_k^{\ast}(z_k) 
\right) y\right]
\nonumber
\\
&=& \sum_{m=0}^{\infty} \int_0^{\infty} e^{-\theta y}
{(\theta y)^m \over m!}dH_k(y) 
\left[ \vc{I} + \theta^{-1}
\left(\vc{C} + \sum_{k \in \K} \vc{D}_k^{\ast}(z_k)\right)
\right]^m 
\nonumber
\\
&=& \sum_{m=0}^{\infty} \gamma_k^{(m)}(\theta)
\left[ \vc{I} + \theta^{-1}
\left(\vc{C} + \sum_{k \in \K} \vc{D}_k^{\ast}(z_k)\right)
\right]^m.
\label{eqn:i-19}
\end{eqnarray}%
\noindent
Substituting (\ref{eqn:a-04}) into (\ref{eqn:i-19}) and 
changing the order of summations, we have
\begin{equation}
\sum_{\svc{n} \in \Z} z_1^{n_1} \cdots z_K^{n_K} 
\vc{A}_k(\vc{n}) 
= \sum_{\svc{n} \in \Z} z_1^{n_1} \cdots z_K^{n_K} 
\sum_{m=0}^{\infty} \gamma_k^{(m)}(\theta) \vc{F}_m(\vc{n}).
\label{eqn:j-18}
\end{equation}%

Comparing the coefficient matrices of $z_1^{n_1} \cdots z_K^{n_K}$ on
both sizes of (\ref{eqn:j-18}), we obtain (\ref{eqn:b-22}).
(\ref{eqn:add-F_0}) is clear from the definition.  The remaining is to
show (\ref{eqn:g-16}). From (\ref{eqn:b-14}) and (\ref{eqn:a-04}), we
have for $m = 0,1,\ldots$,
\begin{eqnarray*}
\lefteqn{\sum_{\svc{n} \in \Z} z_1^{n_1} 
\cdots z_K^{n_K} \vc{F}_{m+1}(\vc{n})} 
\qquad & & 
\\
&=& \sum_{\svc{n} \in \Z} z_1^{n_1} 
\cdots z_K^{n_K} \vc{F}_m(\vc{n})
\left[\vc{I} + \theta^{-1}\left(\vc{C} + 
\sum_{k \in \K} \vc{D}_k^{\ast}(z_k) \right) \right]
\\
&=& \sum_{\svc{n} \in \Z}
z_1^{n_1} \cdots z_K^{n_K} \vc{F}_m(\vc{n}) 
\left(\vc{I} + \theta^{-1}\vc{C} \right) 
\\
& & \quad {} + 
\sum_{\svc{n} \in \Z} \sum_{k \in \K}
\Bigg[ \sum_{l_k=1}^{\infty} 
z_1^{n_1} \cdots z_{k-1}^{n_{k-1}} z_k^{n_k+l_k} 
z_{k+1}^{n_{k+1}} \cdots z_K^{n_K} 
\vc{F}_m(\vc{n}) \cdot \theta^{-1}\vc{D}_k(l_k) \Bigg]
\\
&=& \sum_{\svc{n} \in \Z} 
z_1^{n_1} \cdots z_K^{n_K} \vc{F}_m(\vc{n}) 
\left(\vc{I} + \theta^{-1}\vc{C} \right) 
\\
& & \quad {} +
\sum_{\svc{n} \in \Z^+} z_1^{n_1} \cdots  z_K^{n_K}
\theta^{-1}
\sum_{k \in \K} \sum_{l_k=1}^{n_k} 
\vc{F}_m(\vc{n} - l_k\vc{e}_k)\vc{D}_k(l_k).
\end{eqnarray*}%
\noindent
Comparing the coefficient vectors of $z_1^{n_1} \cdots z_K^{n_K}$ on
both sides of the above equation, we obtain (\ref{eqn:g-16}).
\qed

Next we consider the $\vc{\Gamma}_k(\vc{n})$ in (\ref{eqn:b-13}).
\begin{lem}\label{lemm:a-06}
$\vc{\Gamma}_k(\vc{n})$ ($k \in \K$, $\vc{n} \in \Z$) is determined by
the following recursion:
\begin{eqnarray*}
\vc{\Gamma}_k(\vc{0}) &=& \left[\vc{I} - \vc{P}_k \otimes 
\vc{A}_k(\vc{0}) \right]^{-1}, 
\\
\vc{\Gamma}_k(\vc{n}) 
&=& \sum_{\scriptstyle \svc{0} \le \svc{l} \le \svc{n} 
\atop \scriptstyle \svc{l} \neq \svc{0}} \vc{\Gamma}_k(\vc{n}-\vc{l})
\left[\vc{P}_k \otimes \vc{A}_k(\vc{l}) 
\right]\vc{\Gamma}_k(\vc{0}),
\qquad \vc{n} \in \Z^+.
\end{eqnarray*}%
\end{lem}

\noindent
\textbf{Proof.} 
Note first that (\ref{eqn:b-13}) is equivalent to 
\[
\sum_{\svc{n} \in \Z} z_1^{n_1} \cdots z_K^{n_K} 
\vc{\Gamma}_k(\vc{n}) 
\left[ \vc{I} - \vc{P}_k \otimes \int_0^{\infty} dH_k(y)  
\vc{N}^{\ast}(y, \vc{z}) \right] 
= \vc{I}.
\]%
\noindent
Substituting (\ref{eqn:b-12}) into the above equation, we have
\[
\sum_{\svc{n} \in \Z} z_1^{n_1} \cdots z_K^{n_K} \vc{\Gamma}_k(\vc{n})  
\left[\vc{I} - \vc{P}_k \otimes \sum_{\svc{l} \in \Z}
z_1^{l_1} \cdots z_K^{l_K} 
 \vc{A}_k(\vc{l}) \right] 
= \vc{I},
\]%
\noindent
from which it follows that
\begin{eqnarray*}
\sum_{\svc{n} \in \Z}
z_1^{n_1} \cdots z_K^{n_K} \vc{\Gamma}_k(\vc{n})
- \sum_{\svc{n} \in \Z} 
z_1^{n_1} \cdots z_K^{n_K}
\sum_{\svc{0} \le \svc{l} \le \svc{n}} \vc{\Gamma}_k(\vc{n} - \vc{l})
\left[ \vc{P}_k \otimes \vc{A}_k(\vc{l}) \right] 
= \vc{I}.
\end{eqnarray*}%
\noindent
Comparing the coefficient matrices of $z_1^{n_1} \cdots z_K^{n_K}$ on
both sizes of the above equation, we have
\begin{eqnarray*}
\vc{\Gamma}_k(\vc{0}) - \vc{\Gamma}_k(\vc{0})
\left[ \vc{P}_k \otimes \vc{A}_k(\vc{0}) \right] 
&=& \vc{I}, 
\\
\vc{\Gamma}_k(\vc{n}) - \sum_{\svc{0} \le \svc{l} \le \svc{n}} 
\vc{\Gamma}_k(\vc{n} - \vc{l})
\left[ \vc{P}_k \otimes \vc{A}_k(\vc{l}) \right]
&=& \vc{O},
\qquad \vc{n} \in \Z^{+},
\end{eqnarray*}%
\noindent
or equivalently,
\[
\vc{\Gamma}_k(\vc{0}) 
= \left[\vc{I} - \vc{P}_k 
\otimes \vc{A}_k(\vc{0}) \right]^{-1}, 
\]%
\noindent
and for $\vc{n} \in \Z^{+}$,
\[
\vc{\Gamma}_k(\vc{n}) 
= \sum_{\scriptstyle \svc{0} \le \svc{l} \le \svc{n} 
\atop \scriptstyle \svc{l} \neq \svc{0}} \vc{\Gamma}_k(\vc{n}-\vc{l})
\left[\vc{P}_k \otimes \vc{A}_k(\vc{l}) \right]
\left[\vc{I} - \vc{P}_k 
\otimes \vc{A}_k(\vc{0}) \right]^{-1},
\]%
\noindent
from which Lemma \ref{lemm:a-06} follows. \qed

Finally, we consider the $\vc{v}_k(\vc{n})$ in (\ref{eqn:b-11}). In a
very similar way to derive (\ref{eqn:b-22}), we obtain the following
lemma.
\begin{lem}
$\vc{v}_k(\vc{n})$ ($k \in \K$, $\vc{n} \in \Z$) is given by
\[
\vc{v}_k(\vc{n}) 
= \sum_{m=0}^{\infty} \vc{v}^{(m)}(\theta)
\vc{D}_k\vc{F}_m(\vc{n}),
\]%
\noindent
where $\vc{F}_m(\vc{n})$ is given in (\ref{eqn:add-F_0}) and
(\ref{eqn:g-16}), and
\[
\vc{v}^{(m)}(\theta) 
= \int_0^{\infty} e^{-\theta x}{(\theta x)^m \over m!}
d\vc{v}(x), 
\qquad m=0,1,\dots.
\]%
\end{lem}
Thus $\vc{v}_k(\vc{n})$ is given in terms of the
$\vc{v}^{(m)}(\theta)$.  Because the computation of the
$\vc{v}^{(m)}(\theta)$ has already been studied in \cite{Taki02}, we
summarize the result below. As for the details, readers are referred
to Lemma 3 in \cite{Taki02}.

Note first that 
\begin{equation}
\sum_{m=0}^{\infty} z^m \vc{v}^{(m)}(\theta) 
= \vc{v}^{\ast}(\theta -\theta z),
\label{eqn:f-13}
\end{equation}%
\noindent
where $\vc{v}^{\ast}(s)$ is given in (\ref{eqn:a-11}). Thus,
substituting $\theta - \theta z$ for $s$ in (\ref{eqn:a-11}) and using
(\ref{eqn:f-13}) yield
\begin{eqnarray}
\sum_{m=0}^{\infty} z^m \vc{v}^{(m)}(\theta) 
\left[ (\theta - \theta z) \vc{I} + \vc{C} 
+ \sum_{m=0}^{\infty} z^m  \vc{D}^{(m)}(\theta)  \right]
&=& 
(\theta - \theta z) (1 - \rho)\vc{\kappa},
\label{eqn:j-03}
\end{eqnarray}%
\noindent
where $\vc{D}^{(m)}(\theta)$ denotes
\[
\vc{D}^{(m)}(\theta) 
= \int_0^{\infty}e^{-\theta x} 
{(\theta x)^m \over m!} d\vc{D}(x).
\]%
\noindent
Comparing the coefficient vectors of $z^m$ $(m=0,1,\dots)$ on both
sides of (\ref{eqn:j-03}), we can show that the $\vc{v}^{(m)}(\theta)$
is identical to the steady-state solution of a Markov chain of M/G/1
type whose transition probability matrix is given by \cite{Taki02}
\[
\left[
\begin{array}{ccccc}
\vc{B}_0 + \vc{B}_1 & \vc{B}_2 & \vc{B}_3 & \vc{B}_4 & \cdots
\\
\vc{B}_0            & \vc{B}_1 & \vc{B}_2 & \vc{B}_3 & \cdots
\\
\vc{O}              & \vc{B}_0 & \vc{B}_1 & \vc{B}_2 & \cdots
\\	
\vc{O}              & \vc{O}   & \vc{B}_0 & \vc{B}_1 & \cdots
\\
\vc{O}              & \vc{O}   & \vc{O}   & \vc{B}_0 & \cdots
\\
\vdots              & \vdots   & \vdots   & \vdots   & \ddots
\end{array}
\right],
\]%
\noindent
where
\[
\vc{B}_0 = \vc{I} + \theta^{-1}(\vc{C} + \vc{D}^{(0)}(\theta)),
\quad
\vc{B}_m = \theta^{-1}\vc{D}^{(m)}(\theta), \quad
m \ge 1.
\]%
\noindent
Thus applying the general theory of Markov chains of M/G/1 type
\cite{Neut89}, we can compute the $\vc{v}^{(m)}(\theta)$. As for the
truncation and stopping criteria in computing the steady-state
solution of Markov chains of M/G/1 type, readers are referred to
\cite{Neut89,Taki94b}. 

Let $\vc{d}_k^{(m)}(\theta)$ ($k \in \K$, $m=0,1,\dots$) denote a $1
\times M_k$ vector which satisfies
\begin{eqnarray}
\sum_{m=0}^{\infty} z^m \vc{d}_k^{(m)}(\theta) 
&=& H_k^{\ast}(\theta- \theta z) \vc{\alpha}_k 
(\vc{I} - \vc{P}_k) 
[\vc{I} - H_k^{\ast}(\theta- \theta z) \vc{P}_k]^{-1}.
\label{eqn:i-13}
\end{eqnarray}%
\begin{lem}\label{lem:D}
Under Assumption \ref{ass:a-01}, $\vc{D}^{(m)}(\theta)$ is given by
\[
\vc{D}^{(m)}(\theta) = \sum_{k \in \K} \vc{d}_k^{(m)}(\theta) 
\vc{e} \vc{D}_k,
\qquad m=0,1,\ldots, 
\]%
\noindent
where $\vc{d}_k^{(m)}(\theta)$'s ($k \in \K$) are recursively
determined by
\[
\vc{d}_k^{(0)}(\theta) 
= \gamma_k^{(0)}(\theta) \vc{\alpha}_k 
\left( \vc{I} - \vc{P}_k \right)
\left[ \vc{I} - \gamma_k^{(0)}(\theta) \vc{P}_k \right]^{-1}, 
\]%
\noindent
and for $m =1,2,\dots$,
\begin{eqnarray*}
\vc{d}_k^{(m)}(\theta) 
&=& { \gamma_k^{(m)}(\theta) \over \gamma_k^{(0)}(\theta)}
\vc{d}_k^{(0)}(\theta) 
+
\left[  \sum_{l=1}^m 
\gamma_k^{(l)}(\theta)\vc{d}_k^{(m-l)}(\theta)
\right] \vc{P}_k 
\left[ \vc{I} - \gamma_k^{(0)}(\theta) \vc{P}_k \right]^{-1}.
\end{eqnarray*}%
\end{lem}

\noindent
{\bf Proof}.
Note first that 
\[
\sum_{m=0}^{\infty} z^m \vc{D}^{(m)}(\theta) 
= \vc{D}^{\ast}(\theta - \theta z),
\]%
\noindent
where $\vc{D}^{\ast}(s)$ is given in (\ref{eqn:i-14}). Thus,
substituting $\theta -\theta z$ for $s$ in (\ref{eqn:i-14}) and using
(\ref{eqn:f-14}) and (\ref{eqn:a-12}), we have
\begin{eqnarray}
\sum_{m=0}^{\infty} z^m \vc{D}^{(m)}(\theta) 
&=& \sum_{k \in \K} \sum_{n=1}^{\infty} 
\vc{\alpha}_k \vc{P}_k^{n-1}(\vc{I} - \vc{P}_k)
\vc{e} \{H_k^{\ast}(\theta- \theta z)\}^n \vc{D}_k 
\nonumber
\\
&=& \sum_{k \in \K} H_k^{\ast}(\theta- \theta z) \vc{\alpha}_k 
(\vc{I} - \vc{P}_k) 
[\vc{I} - H_k^{\ast}(\theta- \theta z) \vc{P}_k]^{-1}
\vc{e}\vc{D}_k.
\nonumber
\\
\label{eqn:b-19}
\end{eqnarray}%
\noindent
It then follows from (\ref{eqn:i-13}) and (\ref{eqn:b-19}) that
\begin{eqnarray}
\sum_{m=0}^{\infty} z^m \vc{D}^{(m)}(\theta) 
&=& \sum_{k \in \K} \sum_{m=0}^{\infty} z^m 
\vc{d}_k^{(m)}(\theta)\vc{e}\vc{D}_k
\nonumber
\\
&=& \sum_{m=0}^{\infty} z^m 
\sum_{k \in \K} \vc{d}_k^{(m)}(\theta)\vc{e}\vc{D}_k.
\label{eqn:m-01}
\end{eqnarray}%
\noindent
Note here that 
\begin{equation}
H_k^{\ast}(\theta - \theta z) 
= \sum_{m=0}^{\infty} z^m \gamma_k^{(m)}(\theta),
\qquad k \in \K. 
\label{eqn:c-02}
\end{equation}%
\noindent
Thus from (\ref{eqn:i-13}) and (\ref{eqn:c-02}), we have
\begin{eqnarray*}
\sum_{m=0}^{\infty} z^m \vc{d}_k^{(m)}(\theta) 
\left[ \vc{I} -  \sum_{l=0}^{\infty} z^l 
\gamma_k^{(l)}(\theta) \vc{P}_k \right]
&=& \sum_{m=0}^{\infty} z^m \gamma_k^{(m)}(\theta) 
\vc{\alpha}_k \left[ \vc{I} - \vc{P}_k \right],
\end{eqnarray*}%
\noindent
or equivalently,
\begin{eqnarray*}
\sum_{m=0}^{\infty} z^m \vc{d}_k^{(m)}(\theta)
- \sum_{m=0}^{\infty} z^m  
\sum_{l=0}^m \vc{d}_k^{(m-l)}(\theta)
\gamma_k^{(l)}(\theta)\vc{P}_k
& = &
\sum_{m=0}^{\infty} z^m \gamma_k^{(m)}(\theta) 
\vc{\alpha}_k \left[ \vc{I} - \vc{P}_k \right].
\end{eqnarray*}%
\noindent
Comparing the coefficient vectors of $z^m$ $(m=0,1,\dots)$ on both
sides of the above equation, we have
\begin{equation}
\vc{d}_k^{(0)}(\theta) 
\left[ \vc{I} - \gamma_k^{(0)}(\theta) \vc{P}_k \right]
= \gamma_k^{(0)}(\theta) 
\vc{\alpha}_k \left[ \vc{I} - \vc{P}_k \right], 
\label{eqn:m-02}
\end{equation}%
\noindent
and for $m =1,2,\dots$,
\begin{eqnarray}
\vc{d}_k^{(m)}(\theta)\left[ \vc{I} - \gamma_k^{(0)}(\theta)\vc{P}_k 
\right] 
-
\sum_{l=1}^m \vc{d}_k^{(m-l)}(\theta)
\gamma_k^{(l)}(\theta)\vc{P}_k 
=
\gamma_k^{(m)}(\theta) 
\vc{\alpha}_k \left[ \vc{I} - \vc{P}_k \right].
\label{eqn:m-03}
\end{eqnarray}%
\noindent
Lemma \ref{lem:D} now follows from (\ref{eqn:m-01}), (\ref{eqn:m-02})
and (\ref{eqn:m-03}).  \qed

\section{Implementations of Recursions}\label{sec:5}

In this section, we consider the implementation of recursions for
$\vc{A}_k(\vc{n})$, $\vc{v}_k(\vc{n})$ and $\vc{\Gamma}_k(\vc{n})$,
derived in the preceding section. At a glance, they would seem to be
easy to implement. Contrary to the single arrival case
\cite{Taki01b,Taki02}, however, the computation of the
$\vc{F}_m(\vc{n})$ appeared in $\vc{A}_k(\vc{n})$ and
$\vc{v}_k(\vc{n})$ is not straightforward, because the direct
implementation of the recursion requires very huge memory space and
time-consuming. In what follows, we construct a numerically feasible
procedure to compute the approximate sequences of $\vc{A}_k(\vc{n})$
and $\vc{v}_k(\vc{n})$, avoiding the computation of
$\vc{F}_m(\protect\vc{n})$'s whose contributions to $\vc{A}_k(\vc{n})$
and $\vc{v}_k(\vc{n})$ are negligible, and establish the
truncation/stopping criteria and error bounds.  Further, we propose a
computational procedure for the $\vc{\Gamma}_k(\vc{n})$ and establish
the error bound.

We start with $\vc{A}_k(\vc{n})$ and $\vc{v}_k(\vc{n})$. Note
first that for $k \in \K$,
\[
\sum_{\svc{n} \in \Z} \vc{A}_k(\vc{n}) \vc{e}
= \vc{e}, 
\qquad
\sum_{\svc{n} \in \Z} \vc{v}_k(\vc{n}) \vc{e}
=  \lambda_k^{({\rm B})},
\]%
\noindent
where
\[
\lambda_k^{({\rm B})} = \vc{\pi} \vc{D}_k \vc{e} .
\]%
\noindent
In numerical computation, we have to stop the computation of those
sequences. Thus we develop a numerical procedure to obtain
approximations $\breve{\vc{A}}_k(\vc{n})$ and
$\breve{\vc{v}}_k(\vc{n})$ to $\vc{A}_k(\vc{n})$ and
$\vc{v}_k(\vc{n})$, respectively, while ensuring the following error
bounds: For a given $\varepsilon$ ($0 < \varepsilon < 1$), there exist
$n_A(k)$ and $n_v(k)$ such that
\begin{eqnarray}
\sum_{\scriptstyle \svc{n} \in \Z \atop\scriptstyle 
|\svc{n}| \le n_A(k)}\breve{\vc{A}}_k(\vc{n}) \vc{e}
&>& (1 - \varepsilon)\vc{e},
\label{eqn:e-06}
\\
\sum_{\scriptstyle \svc{n} \in \Z \atop\scriptstyle 
|\svc{n}| \le n_v(k)} \breve{\vc{v}}_k(\vc{n}) \vc{e}
&>& (1 - \varepsilon)\lambda_k^{(\rm{B})},
\label{eqn:e-05}
\end{eqnarray}%
\noindent
where $|\vc{n}| = \sum_{k \in \K} n_k$ for $\vc{n} \in \Z$. In what
follows, we first show our proposed algorithm and then show that the
above error bounds are satisfied.

\subsection*{Numerical algorithm for $\vc{A}_k(\vc{n})$ 
and $\vc{v}_k(\vc{n})$}

\noindent
{\bf Input.}

\medskip

\noindent
Stopping criterion : $\varepsilon$ ($0 < \varepsilon < 1$),

\noindent
Underlying Markov chain 
: $\vc{C}$, $\vc{D}_k$ $(k \in \K)$,

\noindent
Batch size distribution
: $\vc{\alpha}_k$, $\vc{P}_k$ $(k \in \K)$,

\noindent
Service time distribution 
: $H_k(x)$ $(k \in \K)$.

\medskip

\noindent
{\bf Step 1. } 
Choose $\varepsilon_F$ ($0 < \varepsilon_F < 1$) such that
\begin{equation}
{\varepsilon_F \over \varepsilon }
< 
\min_{k \in \K}\min \left( {1 \over \theta h_k},
{
\lambda_k^{(\rm{B})}
\over 
\theta \overline{\vc{v}}^{(1)} \vc{D}_k \vc{e}
} \right),
\label{eqn:j-13}
\end{equation}%
\noindent
where $\overline{\vc{v}}^{(1)} = - \lim_{s \to 0+} 
d\vc{v}^{\ast}(s)/ds$, whose computational procedure can be found in
\cite{Taki94a}.  Then compute the $\gamma_k^{(m)}(\theta)$ and the
$\vc{v}^{(m)}(\theta)$ until they satisfy
\begin{eqnarray}
\sum_{m=0}^{m_{\gamma}(k)}\gamma_k^{(m)}(\theta)(1 - \varepsilon_F)^m
&>& 1 - \varepsilon, \quad k \in \K,
\label{eqn:f-04}
\\
\sum_{m=0}^{m_v(k)} \vc{v}^{(m)}(\theta)\vc{D}_k \vc{e}(1 - \varepsilon_F)^m
&>& (1 - \varepsilon)\lambda_k^{({\rm B})}, \quad k \in \K,
\label{eqn:i-03}
\end{eqnarray}%
\noindent
for some $m_{\gamma}(k)$ and $m_v(k)$, respectively. 
Define $m_{\max}$ as 
\[
m_{\max} = \max_{k \in \K}\max (m_{\gamma}(k), m_v(k)).
\]%

\noindent
{\bf Step 2. }
Choose $\varepsilon_g$ such that $0 < \varepsilon_g < \varepsilon_F$.
Then compute $g_k(n)$ ($n=1,2,\dots$) by (\ref{eqn:a-12}) until
the $g_k(n)$ satisfies
\begin{equation}
\theta^{-1}\sum_{n=1}^{n_g(k)}g_k(n)\vc{D}_k \vc{e}
>
\theta^{-1}\vc{D}_k\vc{e} - {\varepsilon_g \over K} \vc{e},
\quad k \in \K,
\label{eqn:e-15}
\end{equation}%
\noindent
for some $n_g(k)$.

\medskip

\noindent
{\bf Step 3. }
Compute $\breve{\vc{A}}_k(\vc{n})$ and $\breve{\vc{v}}_k(\vc{n})$ by
the following procedure, where the initial values of
$\breve{\vc{A}}_k(\vc{n})$ and $\breve{\vc{v}}_k(\vc{n})$ ($\vc{n} \in
\Z$) are assumed to be $\vc{O}$ and $\vc{0}$, respectively.

\medskip

\begin{narrow}
{\bf Step (3--a). }
Set $\breve{\vc{F}}_0(\vc{0}) = \vc{I}$ and 
$n_F^{(0)} = 0$. Also set
$\breve{\vc{A}}_k(\vc{0}) = \gamma_k^{(0)}(\theta)\vc{I}$ and
$\breve{\vc{v}}_k(\vc{0}) = \vc{v}^{(0)}(\theta)\vc{D}_k$ for all 
$k \in \K$.
\end{narrow}

\begin{narrow}
{\bf Step (3--b). } 
Set $n_F^{(1)} = \max_{k \in \K} n_g(k)$ and $m=1$, 
and compute $\breve{\vc{F}}_1(\vc{n})$'s ($|\vc{n}| \le n_F^{(1)}$) by
\begin{eqnarray}
\breve{\vc{F}}_1(\vc{n}) 
&=&
\left\{
\begin{array}{ll}
\vc{I} + \theta^{-1}\vc{C},
& \quad \mbox{if } \vc{n} = \vc{0},
\\
\theta^{-1} g_k(n_k) \vc{D}_k,
& \quad \mbox{if } \vc{n} \in \Z_k(F_1),\ k \in \K,
\\
\vc{O}, 
& \quad \mbox{otherwise},
\end{array}
\right. 
\label{eqn:e-18}
\end{eqnarray}%
\noindent
where
\[
\Z_k(F_1) = 
\{ \vc{n} ; \vc{n} = n_k\vc{e}_k, n_k=1,2,\dots,n_g(k) \},
\qquad k \in \K.
\]
\end{narrow}

\begin{narrow}
{\bf Step (3--c). } 
For each $k \in \K$, if $m \leq m_{\gamma}(k)$,
add $\gamma_k^{(m)}(\theta) \breve{\vc{F}}_m(\vc{n})$ to
$\breve{\vc{A}}_k(\vc{n})$ for all $\vc{n}$ ($|\vc{n}| \leq
n_F^{(m)}$).  Also, for each $k \in \K$, if $m \leq m_{v}(k)$, add
$\vc{v}^{(m)}(\theta)\vc{D}_k \breve{\vc{F}}_m(\vc{n})$ to
$\breve{\vc{v}}_k(\vc{n})$ for all $\vc{n}$ ($|\vc{n}| \leq
n_F^{(m)}$).
\end{narrow}

\begin{narrow}
{\bf Step (3--d). } 
If $m \geq m_{\max}$, stop computing, and
otherwise, add one to $m$ and go to Step (3--e).
\end{narrow}

\begin{narrow}
{\bf Step (3--e). } 
For each $n=0,1,\dots$, compute
$\breve{\vc{F}}_m(\vc{n})$'s ($|\vc{n}| = n$) by
\begin{eqnarray}
\breve{\vc{F}}_{m}(\vc{n})
&=& U\left(n_F^{(m-1)} - |\vc{n}|\right)
\breve{\vc{F}}_{m-1}(\vc{n})(\vc{I} + \theta^{-1}\vc{C}) 
\nonumber
\\
& & {} + \theta^{-1}\sum_{k \in \K} 
\sum_{l_k=1}^{\min(n_k,\ n_g(k))}
U\left(n_F^{(m-1)} - |\vc{n} - l_k \vc{e}_k| \right)
\nonumber
\\
&& \qquad \qquad {} \cdot 
\breve{\vc{F}}_{m-1}(\vc{n} - l_k\vc{e}_k) g_k(l_k)\vc{D}_k,
\label{eqn:e-23}
\end{eqnarray}%
\noindent
until $\breve{\vc{F}}_m(\vc{n})$'s satisfy $\sum_{|\svc{n}| \le
n^{\ast}} \breve{\vc{F}}_m(\vc{n}) \vc{e} > \left(1 - \varepsilon_F
\right)^m \vc{e}$ for some $n^{\ast}$, where $U(x)$ denotes a unit
step function:
\[
U(x) = \left\{\begin{array}{ll}
1, & x \geq 0, \\
0, & x < 0.
\end{array}\right.
\]%
\noindent
Let $n_F^{(m)} = n^{\ast}$ and go to Step (3--c).
\end{narrow}

\begin{rem}
Note that $\breve{\vc{A}}_k(\vc{n})$ ($|\vc{n}| \le n_A(k)$) and
$\breve{\vc{v}}_k(\vc{n})$ ($|\vc{n}| \le n_v(k)$) obtained by the
above algorithm satisfy
\begin{eqnarray*}
\breve{\vc{A}}_k(\vc{n}) 
&=& 
\sum_{m=0}^{m_{\gamma}(k)} U\left(n_F^{(m)}- |\vc{n}|\right)
\gamma_k^{(m)}(\theta) \breve{\vc{F}}_m(\vc{n}),
\\
\breve{\vc{v}}_k(\vc{n}) 
&=& 
\sum_{m=0}^{m_v(k)} U\left(n_F^{(m)} - |\vc{n}|\right)
\vc{v}^{(m)}(\theta) \vc{D}_k \breve{\vc{F}}_m(\vc{n}),
\end{eqnarray*}%
\noindent
respectively, where 
\begin{eqnarray}
n_A(k) 
&=& \max \left( n_F^{(m)} ; m=0,1,\dots,m_{\gamma}(k) \right),
\label{def_n_A(k)}
\\
n_v(k) 
&=& \max \left( n_F^{(m)} ; m=0,1,\dots,m_v(k) \right).
\nonumber
\end{eqnarray}%
\end{rem}

\begin{rem}
If we are interested only in the $\vc{p}(\vc{n})$ ($|\vc{n}| \leq
N_p$) for some $N_p$, we do not need to compute
$\breve{\vc{F}}(\vc{n})$ for $\vc{n}$ such that $|\vc{n}| >
N_p$. Thus, in this case, $n_g(k)$ is redefined as $\min(n_g(k), N_p)$
and Step (3--e) is replaced by

\begin{narrow}
{\bf Step (3--e'). } 
For each n ($n=0,1,\dots$), compute
$\breve{\vc{F}}_m(\vc{n})$'s ($|\vc{n}| = n$, $\vc{n} \in \Z$) by
(\ref{eqn:e-23}) until $\breve{\vc{F}}_m(\vc{n})$'s satisfy
\[
\sum_{|\svc{n}| \le n^{\ast}} \breve{\vc{F}}_m(\vc{n}) \vc{e} >
\left(1 - \varepsilon_F \right)^m \vc{e},
\]%
\noindent
for some $n^{\ast}$, or $n = N_p$, whichever occurs first. Let
$n_F^{(m)} = \min(n^{\ast}$, $N_p)$ and go to Step (3--c).
\end{narrow}
This procedure can save the computational cost, while maintaining the
accuracy of the results.
\end{rem}

We now provide two lemmas that ensure the above procedure eventually
stops.
\begin{lem}\label{lemma:b-01}
There exist integers $m_{\gamma}(k)$ and $m_v(k)$ satisfying
(\ref{eqn:f-04}) and (\ref{eqn:i-03}), respectively.
\end{lem}

\noindent
\textbf{Proof.}
Substituting $1 - \varepsilon_F$ for $z$ in (\ref{eqn:c-02}), we have
\begin{equation}
\sum_{m=0}^{\infty} \gamma_k^{(m)}(\theta) (1 - \varepsilon_F)^m
= H_k^{\ast}(\theta \varepsilon_F), 
\label{eqn:g-10}
\end{equation}%
\noindent
where $H_k^{\ast}(s)$ denotes the LST of $H_k(x)$. Similarly, from
(\ref{eqn:f-13}), we have
\begin{equation}
\sum_{m=0}^{\infty} 
\vc{v}^{(m)}(\theta)\vc{D}_k \vc{e}(1 - \varepsilon_F)^m
= \vc{v}^{\ast}(\theta \varepsilon_F) \vc{D}_k \vc{e}.
\label{eqn:g-11}
\end{equation}%
\noindent
Note here that 
\begin{eqnarray}
H_k^{\ast}(\theta \varepsilon_F)
&>& 1 - h_k \cdot (\theta \varepsilon_F),
\qquad k \in \K,
\label{eqn:j-11}
\\
\vc{v}^{\ast}(\theta \varepsilon_F) \vc{D}_k \vc{e}
&>& \lambda_k^{(\rm{B})} - \overline{\vc{v}}^{(1)} \vc{D}_k \vc{e} 
\cdot (\theta \varepsilon_F),
\qquad k \in \K,
\label{eqn:j-12}
\end{eqnarray}%
\noindent
because $H_k^{\ast}(s)$ and each element of $\vc{v}^{\ast}(s)$ are
convex functions of $s$. Note also that (\ref{eqn:j-13}) is equivalent
to
\begin{eqnarray}
1 - h_k \theta\varepsilon_F 
&>& 1 - \varepsilon,
\qquad k \in  \K,
\label{eqn:g-12}
\\
\lambda_k^{(\rm{B})} - \overline{\vc{v}}^{(1)} \vc{D}_k \vc{e} 
\theta \varepsilon_F
&>& (1 - \varepsilon)\lambda_k^{({\rm B})},
\qquad k \in \K.
\label{eqn:g-13}
\end{eqnarray}%
\noindent
It then follows from (\ref{eqn:g-10})--(\ref{eqn:g-13}) that
\begin{eqnarray*}
\sum_{m=0}^{\infty} 
\gamma_k^{(m)}(\theta)(1 - \varepsilon_F)^m
&>& 1 - h_k \theta\varepsilon_F 
> 1 - \varepsilon,
\\
\sum_{m=0}^{\infty} 
\vc{v}^{(m)}(\theta)\vc{D}_k \vc{e}(1 - \varepsilon_F)^m
&>&
\lambda_k^{(\rm{B})} - \overline{\vc{v}}^{(1)} \vc{D}_k \vc{e} 
\theta \varepsilon_F
> (1 - \varepsilon)\lambda_k^{({\rm B})},
\end{eqnarray*}%
\noindent
which complete the proof.
\qed

\begin{lem}\label{lemma:b-02}
There exists an integer $n_F^{(m)}$ such that
\begin{equation}
\sum_{\scriptstyle \svc{n} \in \Z \atop\scriptstyle 
|\svc{n}| \le n_F^{(m)}} \breve{\vc{F}}_m(\vc{n}) \vc{e}
> (1 - \varepsilon_F)^m\vc{e},
\qquad \forall m=1,2,\dots,m_{\max}.
\label{eqn:f-03}
\end{equation}%
\end{lem}

\noindent
\textbf{Proof.}
We first consider the case $m=1$. It follows from (\ref{eqn:e-15}) and
(\ref{eqn:e-18}) that
\begin{eqnarray}
\sum_{\scriptstyle \svc{n} \in \Z \atop\scriptstyle 
|\svc{n}| \le n_F^{(1)}}
\breve{\vc{F}}_1(\vc{n}) \vc{e}
&=& \left[ \vc{I} + \theta^{-1}\vc{C} 
+ \theta^{-1}\sum_{k \in \K} 
\sum_{n_k=1}^{n_g(k)} g_k(n_k)\vc{D}_k \right] \vc{e}\nonumber
\\
&>& \vc{e} + \theta^{-1}\vc{C}\vc{e} 
+  \sum_{k \in \K} 
\left[ \theta^{-1}\vc{D}_k\vc{e} - {\varepsilon_g \over K}\vc{e} \right] 
\nonumber
\\
&=& (1 - \varepsilon_g)\vc{e} > (1 - \varepsilon_F)\vc{e},
\label{error_bound_F_1}
\end{eqnarray}%
\noindent
where we use $(\vc{C}+\vc{D})\vc{e} = \vc{0}$ and $\varepsilon_g <
\varepsilon_F$.

Suppose that for some $m$ ($ 1 \le m \le m_{\rm max}-1$), there exists
an integer $n_F^{(m)}$ such that
\begin{equation}
\sum_{\scriptstyle \svc{n} \in \Z \atop\scriptstyle 
|\svc{n}| \le n_F^{(m)}}
\breve{\vc{F}}_m(\vc{n}) \vc{e} > (1 - \varepsilon_F)^m\vc{e}.
\label{error_bound_F_m}
\end{equation}%
\noindent
Using (\ref{eqn:e-18}) and (\ref{eqn:e-23}), 
we have 
\begin{eqnarray}
\breve{\vc{F}}_{m+1}(\vc{n})
&=& U\left(n_F^{(m)} - |\vc{n}|\right)
\breve{\vc{F}}_{m}(\vc{n})\breve{\vc{F}}_1(\vc{0})
\nonumber
\\
& & \qquad {} + \sum_{k \in \K} 
\sum_{l_k=1}^{\min(n_k, n_g(k))}
U\left(n_F^{(m)} - |\vc{n} - l_k \vc{e}_k| \right)
\nonumber
\\
& & {} \qquad \qquad  {} \cdot
\breve{\vc{F}}_m(\vc{n} - l_k\vc{e}_k) \breve{\vc{F}}_1(l_k\vc{e}_k),
\quad \vc{n} \in \Z.
\label{equation:breve{F}_{m+1}(n)}
\end{eqnarray}%
\noindent
It then follows from (\ref{error_bound_F_1}), (\ref{error_bound_F_m}) and
(\ref{equation:breve{F}_{m+1}(n)}) that
\begin{eqnarray}
\sum_{\scriptstyle \svc{n} \in \Z \atop\scriptstyle 
|\svc{n}| \le n_F^{(m)} + n_F^{(1)}}
\breve{\vc{F}}_{m+1}(\vc{n})\vc{e}
&=& 
\sum_{\scriptstyle \svc{n} \in \Z \atop\scriptstyle 
|\svc{n}| \le n_F^{(m)}} \breve{\vc{F}}_m(\vc{n})
\sum_{\scriptstyle \svc{n} \in \Z \atop\scriptstyle 
|\svc{n}| \le n_F^{(1)}} \breve{\vc{F}}_1(\vc{n})\vc{e}
> (1 - \varepsilon_F)^{m+1}\vc{e}.
\label{eqn:e-25}
\end{eqnarray}%
\noindent
Thus we can choose $n_F^{(m+1)}$ in such a way that
\[
n_F^{(m+1)} \le n_F^{(m)} + n_F^{(1)}, \quad m=1,2,\dots,m_{\max}-1,
\]%
\noindent
which completes the proof.
\qed

\begin{thm}\label{thm:b-02}
For $0 < \varepsilon < 1$, the $\breve{\vc{A}}_k(\vc{n})$ and the
$\breve{\vc{v}}_k(\vc{n})$, computed by Step 3, satisfy error bounds
(\ref{eqn:e-06}) and (\ref{eqn:e-05}), respectively.
\end{thm}

\noindent
\textbf{Proof.} 
Using Lemma \ref{lemma:b-01}, Lemma \ref{lemma:b-02} and 
(\ref{def_n_A(k)}), we obtain
\begin{eqnarray*}
\sum_{\scriptstyle \svc{n} \in \Z \atop\scriptstyle 
|\svc{n}| \le n_A(k)}
\breve{\vc{A}}_k(\vc{n}) \vc{e}
&=& \sum_{\scriptstyle \svc{n} \in \Z \atop\scriptstyle 
|\svc{n}| \le n_A(k)}
\sum_{m=0}^{m_{\gamma}(k)}
U\left(n_F^{(m)} - |\vc{n}|\right)
\gamma_k^{(m)}(\theta)
\breve{\vc{F}}_m(\vc{n}) \vc{e}
\\
&=& \sum_{m=0}^{m_{\gamma}(k)}
\gamma_k^{(m)}(\theta)
\sum_{\scriptstyle\svc{n} \in \Z \atop\scriptstyle 
|\svc{n}| \le n_F^{(m)}}
\breve{\vc{F}}_m(\vc{n}) \vc{e}
\\
&>& \sum_{m=0}^{m_{\gamma}(k)}
\gamma_k^{(m)}(\theta)(1 - \varepsilon_F)^m \vc{e}
> (1 - \varepsilon) \vc{e}
, \qquad k \in \K.
\end{eqnarray*}%
\noindent
In the same way, we can obtain (\ref{eqn:e-05}), so that the
proof for the $\breve{\vc{v}}_k(\vc{n})$ is omitted.  \qed

Finally, we consider the $\vc{\Gamma}_k(\vc{n})$. Note here that
\[
\sum_{\svc{n} \in \Z} \vc{\Gamma}_k(\vc{n}) \vc{e}
= \left\{(\vc{I} - \vc{P}_k)^{-1} \vc{e}(M_k) \right\} 
\otimes \vc{e}(M),
\]%
\noindent
where $\vc{e}(m)$ denotes an $m \times 1$ vector whose elements are
all equal to one. Keeping the above equation in mind, we propose to
compute an approximation $\breve{\vc{\Gamma}}_k(\vc{n})$ to
$\vc{\Gamma}_k(\vc{n})$ in the following way.

\medskip

\noindent
{\bf Step 4.} 
For each $k \in \K$, compute $\breve{\vc{\Gamma}}_k(\vc{n})$'s
($|\vc{n}| = n$) for $n = 0,1,\ldots$ by
\begin{eqnarray}
\breve{\vc{\Gamma}}_k(\vc{0}) 
&=& \left[\vc{I} - \vc{P}_k \otimes \breve{\vc{A}}_k(\vc{0}) \right]^{-1},  
\label{step_4_breve_Gamma_0}
\\
\breve{\vc{\Gamma}}_k(\vc{n}) 
&=& \sum_{\scriptstyle \svc{0} \le \svc{l} \le \svc{n} 
\atop \scriptstyle \svc{l} \neq \svc{0}} 
U\left(n_A(k) - |\vc{l}|\right) 
\breve{\vc{\Gamma}}_k(\vc{n}-\vc{l})
\left[\vc{P}_k \otimes \breve{\vc{A}}_k(\vc{l}) 
\right]\breve{\vc{\Gamma}}_k(\vc{0}),
\quad \vc{n} \in \Z^+,
\label{step_4_breve_Gamma_n}
\end{eqnarray}%
\noindent
until $\breve{\vc{\Gamma}}_k(\vc{n})$'s satisfy
\begin{eqnarray}
\sum_{\scriptstyle \svc{n} \in \Z \atop\scriptstyle 
|\svc{n}| \le n_{\Gamma}(k)} \breve{\vc{\Gamma}}_k(\vc{n})\vc{e}
&>& 
\left\{\left(\vc{I} - \vc{P}_k\right)^{-1} \vc{e}(M_k) \right\} 
\otimes \vc{e}(M)
\nonumber
\\
&& \; {} -
\varepsilon
\left\{\left(\vc{I} - \vc{P}_k\right)^{-2} \vc{P}_k
\vc{e}(M_k) \right\}\otimes \vc{e}(M),
\label{stop-A}
\end{eqnarray}%
\noindent
for some integer $n_{\Gamma}(k)$.

\begin{rem}
Let $G_k$ ($k \in \K$) denote a generic random variable representing a
batch size of class $k$. We then have
\[
(\vc{\alpha}_k \otimes \vc{\pi}) 
\sum_{\svc{n} \in \Z} \vc{\Gamma}_k(\vc{n})\vc{e}
= \E[G_k],
\]%
\noindent
and if (\ref{stop-A}) satisfies for some $n_{\Gamma}(k)$,
\[
(\vc{\alpha}_k \otimes \vc{\pi})
\sum_{\scriptstyle \svc{n} \in \Z \atop\scriptstyle 
|\svc{n}| \leq n_{\Gamma}(k)} \breve{\vc{\Gamma}}_k(\vc{n})\vc{e}
> \E[G_k] - {1 \over 2} \E[G_k(G_k-1)] \varepsilon.
\]%
\end{rem}

\begin{lem}
Suppose $\breve{\vc{A}}_k(\vc{n})$ satisfies (\ref{eqn:e-06}).  Then
there exists $n_{\Gamma}(k)$ satisfying (\ref{stop-A}).
\end{lem}

\noindent
{\bf Proof}.  From (\ref{step_4_breve_Gamma_0}) and
(\ref{step_4_breve_Gamma_n}), it can be seen that
$\breve{\vc{A}}_k(\vc{n})$ $(|\vc{n}| \le n_A(k))$ and
$\breve{\vc{\Gamma}}_k(\vc{n})$ $(\vc{n} \in \Z)$ are related by
\[
\sum_{\svc{n} \in \Z} \breve{\vc{\Gamma}}_k(\vc{n})
=
\sum_{m=0}^{\infty}
\left(\vc{P}_k \otimes 
\sum_{\scriptstyle \svc{l} \in \Z 
\atop \scriptstyle |\svc{l}| \le n_A(k)}  
\breve{\vc{A}}_k(\vc{l}) \right)^m.
\]%
\noindent
Post-multiplying both
sides of the above equation by $\vc{e}=\vc{e}(M_k) \otimes \vc{e}(M)$,
we have
\begin{eqnarray*}
\sum_{\svc{n} \in \Z}
\breve{\vc{\Gamma}}_k(\vc{n}) \vc{e} 
&=&
\sum_{m=0}^{\infty}
\left(\vc{P}_k \otimes 
\sum_{\scriptstyle \svc{l} \in \Z 
\atop \scriptstyle |\svc{l}| \le n_A(k)}  
\breve{\vc{A}}_k(\vc{l}) \right)^m 
\cdot [\vc{e}(M_k) \otimes \vc{e}(M)]
\\
&=& \sum_{m=0}^{\infty} [\vc{P}_k^m \vc{e}(M_k) ]
\otimes 
\left[
\left(
\sum_{\scriptstyle \svc{l} \in \Z 
\atop \scriptstyle |\svc{l}| \le n_A(k)}  
\breve{\vc{A}}_k(\vc{l}) \right)^m \vc{e}(M) \right].
\end{eqnarray*}%
\noindent
Further, using (\ref{eqn:e-06}), we obtain
\begin{eqnarray*}
\sum_{\svc{n} \in \Z}
\breve{\vc{\Gamma}}_k(\vc{n}) \vc{e} 
&>& 
\sum_{m=0}^{\infty} (1 - \varepsilon)^m 
[\vc{P}_k^m \vc{e}(M_k) ] \otimes \vc{e}(M) 
\\
&>& \sum_{m=0}^{\infty}
(1 - m \varepsilon) [ \vc{P}_k^m \vc{e}(M_k) ] \otimes \vc{e}(M)
\\
&=&
\left\{ \left(\vc{I} - \vc{P}_k\right)^{-1} \vc{e}(M_k) \right\} 
\otimes \vc{e}(M)
 -
\varepsilon
\left\{ \left(\vc{I} - \vc{P}_k\right)^{-2} \vc{P}_k
\vc{e}(M_k) \right\} \otimes \vc{e}(M),
\end{eqnarray*}%
\noindent
which completes the proof.
\qed

\section{Numerical Examples}\label{sec:6}

In this section, we show some numerical examples for queues with two
arrival streams. Even though the algorithmic analysis has already been
done for the single arrival cases \cite{Taki01b,Taki02}, no numerical
examples were shown there. Thus the numerical result provided below is
the first report in the literature, as for FIFO queues with Markovian
arrival streams having different service time distributions.

In all numerical examples, the counting process of class $k$ ($k=1,2$)
arrivals follows a batch interrupted Poisson process with
geometrically distributed batch size with mean $g$. Namely, the
counting process of class $k$ ($k=1,2$) is characterized by
($\wtilde{\vc{C}}_k$, $\wtilde{\vc{D}}_k(n)$), where
\begin{eqnarray*}
\wtilde{\vc{C}}_k
 &=& 
\left[
\begin{array}{cc}
-2\lambda_k g^{-1}-0.1 & 0.1
\\
0.1 & -0.1
\end{array}
\right],
\\
\wtilde{\vc{D}}_k(n) 
&=& (1-p) p^{n-1}
\left[
\begin{array}{cc}
2\lambda_k g^{-1} & 0
\\
0 & 0
\end{array}
\right], \quad n=1,2,\ldots,
\end{eqnarray*}%
\noindent
where $p=1 - 1/g$. Note that the arrival rate of class $k$ is fixed to
be $\lambda_k$ regardless of the mean batch size $g$.

We now consider three types of the superposition of these two
streams.

\noindent
{[Case P]}
\[
\vc{C} =
\left[
\begin{array}{cc}
-2(\lambda_1 + \lambda_2) g^{-1} - 0.1 & 0.1
\\
0.1 & -0.1
\end{array}
\right],
\]%
\noindent
and for $n=1,2,\ldots$,
\begin{eqnarray*}
\vc{D}_1(n) 
&=& (1-p)p^{n-1}
\left[
\begin{array}{cc}
2\lambda_1 g^{-1}  & 0
\\
0 & 0
\end{array}
\right],
\\
\vc{D}_2(n) 
&=& (1-p)p^{n-1}
\left[
\begin{array}{cc}
2\lambda_2 g^{-1} & 0
\\
0 & 0
\end{array}
\right].
\end{eqnarray*}%

\noindent
{[Case I]}
\[
\vc{C} = \wtilde{\vc{C}}_1 \oplus  \wtilde{\vc{C}}_2, 
\]%
\noindent
and for $n=1,2,\dots$,
\[
\vc{D}_1(n) 
= \wtilde{\vc{D}}_1(n) \otimes  \vc{I}(2), 
\quad
\vc{D}_2(n) 
= \vc{I}(2) \otimes \wtilde{\vc{D}}_2(n),
\]%
\noindent
where $\oplus$ denotes the Kronecker sum, and 

\noindent
{[Case N]}
\[
\vc{C}
=
\left[
\begin{array}{cc}
-2\lambda_1 g^{-1}-0.1 & 0.1
\\
0.1 & -2\lambda_2 g^{-1}-0.1
\end{array}
\right],
\]%
\noindent
and for $n=1,2,\ldots$,
\begin{eqnarray*}
\vc{D}_1(n) 
&=& (1-p)p^{n-1}
\left[
\begin{array}{cc}
2\lambda_1 g^{-1} & 0
\\
0 & 0
\end{array}
\right],
\\
\vc{D}_2(n) 
&=& (1-p)p^{n-1}
\left[
\begin{array}{cc}
0 & 0
\\
0 & 2\lambda_2 g^{-1}
\end{array}
\right].
\end{eqnarray*}%
\noindent
Note that in Case P, two arrival streams are positively correlated, in
Case I, they are independent each other and in Case N, they are
negatively correlated. As for the service time distributions, we
consider two cases, Case GD (class-dependent service times) and Case
GI (i.i.d.\ service times):

\noindent
[Case GD]
\[
H_1 = 1, \mbox{ with prob. } 1,
\quad 
H_2
= 4, \mbox{ with prob. } 1,
\]%

\noindent
[Case GI]
\[
H_k = 
\left\{
\begin{array}{ll}
1, & \mbox{with prob. } 
\lambda_1 / (\lambda_1 + \lambda_2),
\\
4, & \mbox{with prob. } 
\lambda_2 / (\lambda_1 + \lambda_2),
\end{array}
\right.
\quad k=1,2,
\]%
\noindent
where $H_k$ ($k=1,2$) denotes a generic random variable for a service
time of a class $k$ customer. Note that the overall service time
distributions are identical in both cases. We denote the queueing
model with Case $i$ ($i=$ P, I, N) arrivals and Case $j$ ($j=$ GD, GI)
services by Case ($i,j$).

In what follows, we consider two examples, Examples 1 and 2, within
the above settings. In Example 1, we set $\lambda_1 = \lambda_2 =
0.15$, so that $\rho_1 = 0.15$ and $\rho_2 = 0.6$ in Case ($i$, GD)
($i=$ P, I, N), and $\rho_1 = \rho_2 = 0.375$ in Case ($i$,\ GI) ($i=$
P, I, N).  On the other hand, in Example 2, we set $\lambda_1 = 0.4$
and $\lambda_2 = 0.1$, so that $\rho_1 = \rho_2 = 0.4$ in Case ($i$,
GD) ($i=$ P, I, N) and that $\rho_1 = 0.64$ and $\rho_2 = 0.16$ in
Case ($i$, GI) ($i=$ P, I, N).

\subsection{Efficiency of the algorithm}

Before showing the quantitative behavior of the queue length
distribution, we discuss the efficiency of our numerical algorithm for
the $\vc{F}_m(\vc{n})$. It follows from (\ref{eqn:a-04}) that for $m
=0,1,\dots$,
\begin{equation}
\sum_{\svc{n} \in \Z} \vc{F}_m(\vc{n})
=
\left[\vc{I} + \theta^{-1}
\left(\vc{C} + \sum_{k \in \K} 
\sum_{n_k=1}^{\infty} \vc{D}_k(n_k)\right)
\right]^m, 
\label{add-F_m(n)}
\end{equation}%
\noindent
where $\vc{I} + \theta^{-1}[\vc{C} + \sum_{k \in \K}
\sum_{n_k=1}^{\infty} \vc{D}_k(n_k)]$ is a stochastic matrix. Thus a
straightforward implementation of the recursion for the
$\vc{F}_m(\vc{n})$ in (\ref{eqn:add-F_0}) and (\ref{eqn:g-16}) would
be the following. We first truncate the $\vc{D}_k(n_k)$ at $n_k =
n_g'(k)$ in such a way that
\[
\theta^{-1} \sum_{n_k=1}^{n_g'(k)} \vc{D}_k(n_k) \vc{e}
> 
\theta^{-1} \vc{D}_k \vc{e} 
-
{\varepsilon_g' \over K} \vc{e},
\]%
\noindent
so that
\[
\left[\vc{I} + \theta^{-1}
\left(\vc{C} + \sum_{k \in \K} 
\sum_{n_k=1}^{n_g'(k)} \vc{D}_k(n_k)\right)
\right]\vc{e} > (1-\varepsilon_g') \vc{e}.
\]%
\noindent
We then compute all terms obtained by expanding the right hand side of
(\ref{add-F_m(n)}) with the truncated $\vc{D}_k(n_k)$ ($k \in
\K$). Note that if $\varepsilon_g' = \varepsilon_F$, the resulting
$\vc{F}_m(\vc{n})$ satisfies (\ref{eqn:f-03}) in Lemma
\ref{lemma:b-02}, where the summation on the left hand side of 
(\ref{eqn:f-03}) is taken for all computed $\vc{F}_m(\vc{n})$'s.

In Table \ref{table-Fm}, we show the numbers of $\vc{F}_m(\vc{n})$'s
computed by our algorithm and the above straightforward
implementation, using Example 1, where we set $\varepsilon=10^{-6}$,
$\varepsilon_F= \mbox{r.h.s.\ of (\ref{eqn:j-13})} \times
\varepsilon/2$ and $\varepsilon_g = \varepsilon_F/10$. We observe that
for unbounded batch size cases (i.e., $g > 1$), the number of the
computed $\vc{F}_m(\vc{n})$'s in our algorithm is less than that in
the straightforward algorithm about by three order of magnitude.
Thus, compared to the straightforward implementation, our algorithm is
very efficient in terms of the computational time when the batch size
is unbounded.
\begin{table}[tb]
\begin{center}
\caption{Number of computed $\protect\vc{F}_m(\protect\vc{n})$'s 
in Example 1.}
\label{table-Fm}
A : Our algorithm \quad B : Straightforward
\begin{tabular}{@{}llr@{}lr@{}lr@{}lr@{}l}
\hline 
\multicolumn{1}{c}{Case}
& \multicolumn{1}{c}{} 
& \multicolumn{2}{c}{$g=1$} 
& \multicolumn{2}{c}{$g=2$}
& \multicolumn{2}{c}{$g=5$} 
& \multicolumn{2}{c}{$g=10$}
\\
\hline
\down{2.5mm}{(P, GD)} & 
A 
& 1.021 
& $\times 10^6$ 
& 3.918 
& $ \times 10^6$
& 2.504 
& $\times 10^7$
& 1.257 
& $\times 10^8$
\\
\cline{2-10}
&
B 
& 2.453 
& $\times 10^6$
& 2.107 
& $\times 10^9$
& 3.696 
& $\times 10^{10}$
& 4.243 
& $\times 10^{11}$
\\
\hline
\hline

\down{2.5mm}{(P, GI)}
& A 
& 1.021 
& $\times 10^6$
& 2.898 
& $\times 10^6$
& 1.476 
& $\times 10^7$
& 6.732 
& $\times 10^7$
\\
\cline{2-10}

&
B 
& 2.453 
& $\times 10^6$
& 1.438 
& $\times 10^9$
& 1.812 
& $\times 10^{10}$
& 1.825 
& $\times 10^{11}$
\\
\hline
\hline

\down{2.5mm}{(I, GD)}
&
A 
& 6.108 
& $\times 10^5$
& 3.314 
& $\times 10^6$
& 3.012 
& $\times 10^7$
& 1.854 
& $\times 10^8$
\\
\cline{2-10}
 
&
B 
& 1.993 
& $\times 10^6$
& 2.919 
& $\times 10^9$
& 9.286 
& $\times 10^{10}$
& 1.649 
& $\times 10^{12}$
\\
\hline
\hline
 
\down{2.5mm}{(I, GI)}
&
A 
& 4.123 
& $\times 10^5$
& 1.859 
& $\times 10^6$
& 1.532 
& $\times 10^7$
& 9.032 
& $\times 10^7$
\\
\cline{2-10}

&
B 
& 1.253 
& $\times 10^6$
& 1.292 
& $\times 10^9$
& 3.741 
& $\times 10^{10}$	
& 6.173 
& $\times 10^{11}$	
\\
\hline
\hline

\down{2.5mm}{(N, GD)}
&
A 
& 6.657 
& $\times 10^4$
& 8.895 
& $\times 10^5$
& 1.165 
& $\times 10^7$
& 8.095 
& $\times 10^7$
\\
\cline{2-10}

&
B 
& 1.378 
& $\times 10^5$
& 2.620 
& $\times 10^8$
& 1.066 
& $\times 10^{10}$
& 1.955 
& $\times 10^{11}$
\\
\hline
\hline

\down{2.5mm}{(N, GI)}
&
A 
& 1.411 
& $\times 10^4$
& 3.113 
& $\times 10^5$
& 4.813 
& $\times 10^6$
& 3.540 
& $\times 10^7$
\\
\cline{2-10}
 
&
B 
& 2.743 
& $\times 10^4$
& 7.158 
& $\times 10^7$
& 3.278 
& $\times 10^9$
& 6.442 
& $\times 10^{10}$
\\
\hline
\end{tabular}
\end{center}
\end{table}

We note that a very huge memory space is required to store all
$\breve{\vc{F}}_m(\vc{n})$'s in some cases, even using our truncation
and stopping criteria. For example, in Case (I, GD) with $g=10$, the
memory space to store all $\breve{\vc{F}}_m(\vc{n})$'s is given by $16
\times 1.854 \times 10^8 \times 8$ bytes ${} \approx 23.73$ Gbytes,
because each $\breve{\vc{F}}_m(\vc{n})$ is a $4 \times 4$ matrix and
one element requires 8 bytes in double precision. Thus in our
implementation, every time $\breve{\vc{F}}_m(\vc{n})$'s for each $m$
are obtained, we compute the contributions of
$\breve{\vc{F}}_m(\vc{n})$'s to $\breve{\vc{A}}_k(\vc{n})$ and
$\breve{\vc{v}}_k(\vc{n})$ in Step (3--c), and discard all
$\breve{\vc{F}}_{m-1}(\vc{n})$'s.

Table \ref{table-stored} shows the maximum number of
$\breve{\vc{F}}_m(\vc{n})$'s stored temporarily in our algorithm,
where the ratio of it to the total number of computed
$\breve{\vc{F}}_m(\vc{n})$'s is also shown in parenthesis. We observe
that in most cases, the number of temporarily stored
$\breve{\vc{F}}_m(\vc{n})$'s is a few percent of the total number of
computed ones. Thus our implementation is expected to save the
required memory space, especially when a large number of
$\breve{\vc{F}}_m(\vc{n})$'s should be computed.
\begin{table}[tb]
\begin{center}
\caption{Number of 
stored $\breve{\protect\vc{F}}_m(\protect\vc{n})$'s in Example 1.}
\label{table-stored}
\begin{tabular}{@{}ccccc}
\hline 
Case
& $g=1$ 
& $g=2$
& $g=5$ 
& $g=10$
\\
\hline
\down{2.5mm}{(P, GD)}

& 27225 
& 90601
& 455625
& 1651227 
\\
& (2.67\%)
& (2.31\%)
& (1.82\%)
& (1.31\%)
\\
\hline
\down{2.5mm}{(P, GI)}
& 27225
& 75076
& 330051
& 1125723
\\
& (2.67\%)
& (2.59\%)
& (2.24\%)
& (1.67\%)
\\
\hline
\down{2.5mm}{(I, GD)}
& 16641
& 68121
& 399424
& 1548781
\\
& (2.72\%)
& (2.06\%)
& (1.33\%)
& (0.84\%)
\\
\hline
\down{2.5mm}{(I, GI)}
& 13110 
& 47524
& 263683 
& \phantom{1}997003 
\\
& (3.18\%)
& (2.56\%)
& (1.72\%)
& (1.10\%)
\\
\hline
\down{2.5mm}{(N, GD)}
& \phantom{1}4970 
& 41209
& 324331 
& 1387686
\\
& (7.47\%)
& (4.63\%)
& (2.78\%)
& (1.71\%)
\\
\hline
\down{2.5mm}{(N, GI)}
& \phantom{1}1764 
& 21171
& 187491 
& \phantom{1}833571
\\
& (12.51\%)\phantom{1}
& (6.80\%)
& (3.90\%)
& (2.35\%)
\\
\hline
\end{tabular}
\end{center}
\end{table}

\subsection{Number of customers in Example 1}

Figures \ref{figure-A}--\ref{figure-C} plot the complementary
distributions of the total number $N$ of customers in Case ($i$, GD)
and Case ($i$, GI) ($i=$ P, I, N), where the batch size is fixed to be
one, i.e., $g=1$. Note that the overall input processes in Case (P,
GD) and Case (P, GI) are identical, so that the distributions of the
total number of customers are also identical, as shown in Figure
\ref{figure-A}.  However, as shown in Table \ref{table-14}, the joint
queue length distributions in these two cases are different. Note also
that in Case (P, GI), $\vc{p}(n_1,n_2)\vc{e} = \vc{p}(n_2,n_1)\vc{e}$,
because the conditional joint distribution $\Pr(N_1=n_1, N_2=n_2 \mid
N_1+N_2 = n_1+n_2)$ follows a binomial distribution with parameter
0.5. We also observe that $\vc{p}(n,n)\vc{e}$'s in both cases take the
same value for each $n$. Unfortunately, we cannot provide any
intuitive explanation of this phenomenon.
\begin{figure}[tb]
\begin{center}
\input{fig_01.tex}
\end{center}
\caption{Complementary distribution of total number of customers 
in Example 1.}
\label{figure-A}
\end{figure}

\begin{figure}[tb]
\begin{center}
% GNUPLOT: LaTeX picture
\setlength{\unitlength}{0.240900pt}
\ifx\plotpoint\undefined\newsavebox{\plotpoint}\fi
\sbox{\plotpoint}{\rule[-0.200pt]{0.400pt}{0.400pt}}%
\begin{picture}(1200,675)(0,0)
\font\gnuplot=cmr10 at 10pt
\gnuplot
\sbox{\plotpoint}{\rule[-0.200pt]{0.400pt}{0.400pt}}%
\put(161.0,123.0){\rule[-0.200pt]{4.818pt}{0.400pt}}
\put(141,123){\makebox(0,0)[r]{0}}
\put(1119.0,123.0){\rule[-0.200pt]{4.818pt}{0.400pt}}
\put(161.0,225.0){\rule[-0.200pt]{4.818pt}{0.400pt}}
\put(141,225){\makebox(0,0)[r]{0.2}}
\put(1119.0,225.0){\rule[-0.200pt]{4.818pt}{0.400pt}}
\put(161.0,328.0){\rule[-0.200pt]{4.818pt}{0.400pt}}
\put(141,328){\makebox(0,0)[r]{0.4}}
\put(1119.0,328.0){\rule[-0.200pt]{4.818pt}{0.400pt}}
\put(161.0,430.0){\rule[-0.200pt]{4.818pt}{0.400pt}}
\put(141,430){\makebox(0,0)[r]{0.6}}
\put(1119.0,430.0){\rule[-0.200pt]{4.818pt}{0.400pt}}
\put(161.0,533.0){\rule[-0.200pt]{4.818pt}{0.400pt}}
\put(141,533){\makebox(0,0)[r]{0.8}}
\put(1119.0,533.0){\rule[-0.200pt]{4.818pt}{0.400pt}}
\put(161.0,635.0){\rule[-0.200pt]{4.818pt}{0.400pt}}
\put(141,635){\makebox(0,0)[r]{1}}
\put(1119.0,635.0){\rule[-0.200pt]{4.818pt}{0.400pt}}
\put(161.0,123.0){\rule[-0.200pt]{0.400pt}{4.818pt}}
\put(161,82){\makebox(0,0){0}}
\put(161.0,615.0){\rule[-0.200pt]{0.400pt}{4.818pt}}
\put(324.0,123.0){\rule[-0.200pt]{0.400pt}{4.818pt}}
\put(324,82){\makebox(0,0){5}}
\put(324.0,615.0){\rule[-0.200pt]{0.400pt}{4.818pt}}
\put(487.0,123.0){\rule[-0.200pt]{0.400pt}{4.818pt}}
\put(487,82){\makebox(0,0){10}}
\put(487.0,615.0){\rule[-0.200pt]{0.400pt}{4.818pt}}
\put(650.0,123.0){\rule[-0.200pt]{0.400pt}{4.818pt}}
\put(650,82){\makebox(0,0){15}}
\put(650.0,615.0){\rule[-0.200pt]{0.400pt}{4.818pt}}
\put(813.0,123.0){\rule[-0.200pt]{0.400pt}{4.818pt}}
\put(813,82){\makebox(0,0){20}}
\put(813.0,615.0){\rule[-0.200pt]{0.400pt}{4.818pt}}
\put(976.0,123.0){\rule[-0.200pt]{0.400pt}{4.818pt}}
\put(976,82){\makebox(0,0){25}}
\put(976.0,615.0){\rule[-0.200pt]{0.400pt}{4.818pt}}
\put(1139.0,123.0){\rule[-0.200pt]{0.400pt}{4.818pt}}
\put(1139,82){\makebox(0,0){30}}
\put(1139.0,615.0){\rule[-0.200pt]{0.400pt}{4.818pt}}
\put(161.0,123.0){\rule[-0.200pt]{235.600pt}{0.400pt}}
\put(1139.0,123.0){\rule[-0.200pt]{0.400pt}{123.341pt}}
\put(161.0,635.0){\rule[-0.200pt]{235.600pt}{0.400pt}}
\put(40,379){\makebox(0,0){Prob.}}
\put(650,21){\makebox(0,0){Total number of customers}}
\put(161.0,123.0){\rule[-0.200pt]{0.400pt}{123.341pt}}
\put(979,585){\makebox(0,0)[r]{Case (I, GD)}}
\put(999.0,585.0){\rule[-0.200pt]{24.090pt}{0.400pt}}
\put(161,507){\usebox{\plotpoint}}
\multiput(161.58,503.21)(0.497,-1.020){63}{\rule{0.120pt}{0.912pt}}
\multiput(160.17,505.11)(33.000,-65.107){2}{\rule{0.400pt}{0.456pt}}
\multiput(194.58,436.83)(0.497,-0.831){61}{\rule{0.120pt}{0.763pt}}
\multiput(193.17,438.42)(32.000,-51.417){2}{\rule{0.400pt}{0.381pt}}
\multiput(226.58,384.42)(0.497,-0.652){63}{\rule{0.120pt}{0.621pt}}
\multiput(225.17,385.71)(33.000,-41.711){2}{\rule{0.400pt}{0.311pt}}
\multiput(259.58,341.72)(0.497,-0.562){61}{\rule{0.120pt}{0.550pt}}
\multiput(258.17,342.86)(32.000,-34.858){2}{\rule{0.400pt}{0.275pt}}
\multiput(291.00,306.92)(0.549,-0.497){57}{\rule{0.540pt}{0.120pt}}
\multiput(291.00,307.17)(31.879,-30.000){2}{\rule{0.270pt}{0.400pt}}
\multiput(324.00,276.92)(0.661,-0.497){47}{\rule{0.628pt}{0.120pt}}
\multiput(324.00,277.17)(31.697,-25.000){2}{\rule{0.314pt}{0.400pt}}
\multiput(357.00,251.92)(0.765,-0.496){39}{\rule{0.710pt}{0.119pt}}
\multiput(357.00,252.17)(30.527,-21.000){2}{\rule{0.355pt}{0.400pt}}
\multiput(389.00,230.92)(0.923,-0.495){33}{\rule{0.833pt}{0.119pt}}
\multiput(389.00,231.17)(31.270,-18.000){2}{\rule{0.417pt}{0.400pt}}
\multiput(422.00,212.92)(1.079,-0.494){27}{\rule{0.953pt}{0.119pt}}
\multiput(422.00,213.17)(30.021,-15.000){2}{\rule{0.477pt}{0.400pt}}
\multiput(454.00,197.92)(1.401,-0.492){21}{\rule{1.200pt}{0.119pt}}
\multiput(454.00,198.17)(30.509,-12.000){2}{\rule{0.600pt}{0.400pt}}
\multiput(487.00,185.92)(1.694,-0.491){17}{\rule{1.420pt}{0.118pt}}
\multiput(487.00,186.17)(30.053,-10.000){2}{\rule{0.710pt}{0.400pt}}
\multiput(520.00,175.93)(1.834,-0.489){15}{\rule{1.522pt}{0.118pt}}
\multiput(520.00,176.17)(28.841,-9.000){2}{\rule{0.761pt}{0.400pt}}
\multiput(552.00,166.93)(2.476,-0.485){11}{\rule{1.986pt}{0.117pt}}
\multiput(552.00,167.17)(28.879,-7.000){2}{\rule{0.993pt}{0.400pt}}
\multiput(585.00,159.93)(2.399,-0.485){11}{\rule{1.929pt}{0.117pt}}
\multiput(585.00,160.17)(27.997,-7.000){2}{\rule{0.964pt}{0.400pt}}
\multiput(617.00,152.93)(3.604,-0.477){7}{\rule{2.740pt}{0.115pt}}
\multiput(617.00,153.17)(27.313,-5.000){2}{\rule{1.370pt}{0.400pt}}
\multiput(650.00,147.94)(4.722,-0.468){5}{\rule{3.400pt}{0.113pt}}
\multiput(650.00,148.17)(25.943,-4.000){2}{\rule{1.700pt}{0.400pt}}
\multiput(683.00,143.94)(4.575,-0.468){5}{\rule{3.300pt}{0.113pt}}
\multiput(683.00,144.17)(25.151,-4.000){2}{\rule{1.650pt}{0.400pt}}
\multiput(715.00,139.95)(7.160,-0.447){3}{\rule{4.500pt}{0.108pt}}
\multiput(715.00,140.17)(23.660,-3.000){2}{\rule{2.250pt}{0.400pt}}
\put(748,136.17){\rule{6.500pt}{0.400pt}}
\multiput(748.00,137.17)(18.509,-2.000){2}{\rule{3.250pt}{0.400pt}}
\put(780,134.17){\rule{6.700pt}{0.400pt}}
\multiput(780.00,135.17)(19.094,-2.000){2}{\rule{3.350pt}{0.400pt}}
\put(813,132.17){\rule{6.700pt}{0.400pt}}
\multiput(813.00,133.17)(19.094,-2.000){2}{\rule{3.350pt}{0.400pt}}
\put(846,130.67){\rule{7.709pt}{0.400pt}}
\multiput(846.00,131.17)(16.000,-1.000){2}{\rule{3.854pt}{0.400pt}}
\put(878,129.17){\rule{6.700pt}{0.400pt}}
\multiput(878.00,130.17)(19.094,-2.000){2}{\rule{3.350pt}{0.400pt}}
\put(911,127.67){\rule{7.709pt}{0.400pt}}
\multiput(911.00,128.17)(16.000,-1.000){2}{\rule{3.854pt}{0.400pt}}
\put(943,126.67){\rule{7.950pt}{0.400pt}}
\multiput(943.00,127.17)(16.500,-1.000){2}{\rule{3.975pt}{0.400pt}}
\put(1009,125.67){\rule{7.709pt}{0.400pt}}
\multiput(1009.00,126.17)(16.000,-1.000){2}{\rule{3.854pt}{0.400pt}}
\put(976.0,127.0){\rule[-0.200pt]{7.950pt}{0.400pt}}
\put(1074,124.67){\rule{7.709pt}{0.400pt}}
\multiput(1074.00,125.17)(16.000,-1.000){2}{\rule{3.854pt}{0.400pt}}
\put(1041.0,126.0){\rule[-0.200pt]{7.950pt}{0.400pt}}
\put(1106.0,125.0){\rule[-0.200pt]{7.950pt}{0.400pt}}
\sbox{\plotpoint}{\rule[-0.400pt]{0.800pt}{0.800pt}}%
\put(979,524){\makebox(0,0)[r]{Case (I, GI)}}
\put(999.0,524.0){\rule[-0.400pt]{24.090pt}{0.800pt}}
\put(161,507){\usebox{\plotpoint}}
\multiput(162.41,498.62)(0.503,-1.149){59}{\rule{0.121pt}{2.018pt}}
\multiput(159.34,502.81)(33.000,-70.811){2}{\rule{0.800pt}{1.009pt}}
\multiput(195.41,425.05)(0.503,-0.929){57}{\rule{0.121pt}{1.675pt}}
\multiput(192.34,428.52)(32.000,-55.523){2}{\rule{0.800pt}{0.838pt}}
\multiput(227.41,367.44)(0.503,-0.714){59}{\rule{0.121pt}{1.339pt}}
\multiput(224.34,370.22)(33.000,-44.220){2}{\rule{0.800pt}{0.670pt}}
\multiput(260.41,321.23)(0.503,-0.593){57}{\rule{0.121pt}{1.150pt}}
\multiput(257.34,323.61)(32.000,-35.613){2}{\rule{0.800pt}{0.575pt}}
\multiput(291.00,286.09)(0.548,-0.503){53}{\rule{1.080pt}{0.121pt}}
\multiput(291.00,286.34)(30.758,-30.000){2}{\rule{0.540pt}{0.800pt}}
\multiput(324.00,256.09)(0.662,-0.504){43}{\rule{1.256pt}{0.121pt}}
\multiput(324.00,256.34)(30.393,-25.000){2}{\rule{0.628pt}{0.800pt}}
\multiput(357.00,231.09)(0.808,-0.505){33}{\rule{1.480pt}{0.122pt}}
\multiput(357.00,231.34)(28.928,-20.000){2}{\rule{0.740pt}{0.800pt}}
\multiput(389.00,211.09)(1.055,-0.507){25}{\rule{1.850pt}{0.122pt}}
\multiput(389.00,211.34)(29.160,-16.000){2}{\rule{0.925pt}{0.800pt}}
\multiput(422.00,195.09)(1.179,-0.509){21}{\rule{2.029pt}{0.123pt}}
\multiput(422.00,195.34)(27.790,-14.000){2}{\rule{1.014pt}{0.800pt}}
\multiput(454.00,181.08)(1.586,-0.512){15}{\rule{2.600pt}{0.123pt}}
\multiput(454.00,181.34)(27.604,-11.000){2}{\rule{1.300pt}{0.800pt}}
\multiput(487.00,170.08)(1.999,-0.516){11}{\rule{3.133pt}{0.124pt}}
\multiput(487.00,170.34)(26.497,-9.000){2}{\rule{1.567pt}{0.800pt}}
\multiput(520.00,161.08)(2.665,-0.526){7}{\rule{3.857pt}{0.127pt}}
\multiput(520.00,161.34)(23.994,-7.000){2}{\rule{1.929pt}{0.800pt}}
\multiput(552.00,154.07)(3.476,-0.536){5}{\rule{4.600pt}{0.129pt}}
\multiput(552.00,154.34)(23.452,-6.000){2}{\rule{2.300pt}{0.800pt}}
\multiput(585.00,148.06)(4.958,-0.560){3}{\rule{5.320pt}{0.135pt}}
\multiput(585.00,148.34)(20.958,-5.000){2}{\rule{2.660pt}{0.800pt}}
\put(617,141.34){\rule{6.800pt}{0.800pt}}
\multiput(617.00,143.34)(18.886,-4.000){2}{\rule{3.400pt}{0.800pt}}
\put(650,137.84){\rule{7.950pt}{0.800pt}}
\multiput(650.00,139.34)(16.500,-3.000){2}{\rule{3.975pt}{0.800pt}}
\put(683,134.84){\rule{7.709pt}{0.800pt}}
\multiput(683.00,136.34)(16.000,-3.000){2}{\rule{3.854pt}{0.800pt}}
\put(715,132.34){\rule{7.950pt}{0.800pt}}
\multiput(715.00,133.34)(16.500,-2.000){2}{\rule{3.975pt}{0.800pt}}
\put(748,130.34){\rule{7.709pt}{0.800pt}}
\multiput(748.00,131.34)(16.000,-2.000){2}{\rule{3.854pt}{0.800pt}}
\put(780,128.84){\rule{7.950pt}{0.800pt}}
\multiput(780.00,129.34)(16.500,-1.000){2}{\rule{3.975pt}{0.800pt}}
\put(813,127.34){\rule{7.950pt}{0.800pt}}
\multiput(813.00,128.34)(16.500,-2.000){2}{\rule{3.975pt}{0.800pt}}
\put(846,125.84){\rule{7.709pt}{0.800pt}}
\multiput(846.00,126.34)(16.000,-1.000){2}{\rule{3.854pt}{0.800pt}}
\put(911,124.84){\rule{7.709pt}{0.800pt}}
\multiput(911.00,125.34)(16.000,-1.000){2}{\rule{3.854pt}{0.800pt}}
\put(943,123.84){\rule{7.950pt}{0.800pt}}
\multiput(943.00,124.34)(16.500,-1.000){2}{\rule{3.975pt}{0.800pt}}
\put(878.0,127.0){\rule[-0.400pt]{7.950pt}{0.800pt}}
\put(1041,122.84){\rule{7.950pt}{0.800pt}}
\multiput(1041.00,123.34)(16.500,-1.000){2}{\rule{3.975pt}{0.800pt}}
\put(976.0,125.0){\rule[-0.400pt]{15.658pt}{0.800pt}}
\put(1074.0,124.0){\rule[-0.400pt]{15.658pt}{0.800pt}}
\end{picture}
\end{center}
\caption{Complementary distribution of total number of customers 
in Example 1.}
\label{figure-B}
\end{figure}

\begin{figure}[tb]
\begin{center}
% GNUPLOT: LaTeX picture
\setlength{\unitlength}{0.240900pt}
\ifx\plotpoint\undefined\newsavebox{\plotpoint}\fi
\sbox{\plotpoint}{\rule[-0.200pt]{0.400pt}{0.400pt}}%
\begin{picture}(1200,675)(0,0)
\font\gnuplot=cmr10 at 10pt
\gnuplot
\sbox{\plotpoint}{\rule[-0.200pt]{0.400pt}{0.400pt}}%
\put(161.0,123.0){\rule[-0.200pt]{4.818pt}{0.400pt}}
\put(141,123){\makebox(0,0)[r]{0}}
\put(1119.0,123.0){\rule[-0.200pt]{4.818pt}{0.400pt}}
\put(161.0,225.0){\rule[-0.200pt]{4.818pt}{0.400pt}}
\put(141,225){\makebox(0,0)[r]{0.2}}
\put(1119.0,225.0){\rule[-0.200pt]{4.818pt}{0.400pt}}
\put(161.0,328.0){\rule[-0.200pt]{4.818pt}{0.400pt}}
\put(141,328){\makebox(0,0)[r]{0.4}}
\put(1119.0,328.0){\rule[-0.200pt]{4.818pt}{0.400pt}}
\put(161.0,430.0){\rule[-0.200pt]{4.818pt}{0.400pt}}
\put(141,430){\makebox(0,0)[r]{0.6}}
\put(1119.0,430.0){\rule[-0.200pt]{4.818pt}{0.400pt}}
\put(161.0,533.0){\rule[-0.200pt]{4.818pt}{0.400pt}}
\put(141,533){\makebox(0,0)[r]{0.8}}
\put(1119.0,533.0){\rule[-0.200pt]{4.818pt}{0.400pt}}
\put(161.0,635.0){\rule[-0.200pt]{4.818pt}{0.400pt}}
\put(141,635){\makebox(0,0)[r]{1}}
\put(1119.0,635.0){\rule[-0.200pt]{4.818pt}{0.400pt}}
\put(161.0,123.0){\rule[-0.200pt]{0.400pt}{4.818pt}}
\put(161,82){\makebox(0,0){0}}
\put(161.0,615.0){\rule[-0.200pt]{0.400pt}{4.818pt}}
\put(324.0,123.0){\rule[-0.200pt]{0.400pt}{4.818pt}}
\put(324,82){\makebox(0,0){5}}
\put(324.0,615.0){\rule[-0.200pt]{0.400pt}{4.818pt}}
\put(487.0,123.0){\rule[-0.200pt]{0.400pt}{4.818pt}}
\put(487,82){\makebox(0,0){10}}
\put(487.0,615.0){\rule[-0.200pt]{0.400pt}{4.818pt}}
\put(650.0,123.0){\rule[-0.200pt]{0.400pt}{4.818pt}}
\put(650,82){\makebox(0,0){15}}
\put(650.0,615.0){\rule[-0.200pt]{0.400pt}{4.818pt}}
\put(813.0,123.0){\rule[-0.200pt]{0.400pt}{4.818pt}}
\put(813,82){\makebox(0,0){20}}
\put(813.0,615.0){\rule[-0.200pt]{0.400pt}{4.818pt}}
\put(976.0,123.0){\rule[-0.200pt]{0.400pt}{4.818pt}}
\put(976,82){\makebox(0,0){25}}
\put(976.0,615.0){\rule[-0.200pt]{0.400pt}{4.818pt}}
\put(1139.0,123.0){\rule[-0.200pt]{0.400pt}{4.818pt}}
\put(1139,82){\makebox(0,0){30}}
\put(1139.0,615.0){\rule[-0.200pt]{0.400pt}{4.818pt}}
\put(161.0,123.0){\rule[-0.200pt]{235.600pt}{0.400pt}}
\put(1139.0,123.0){\rule[-0.200pt]{0.400pt}{123.341pt}}
\put(161.0,635.0){\rule[-0.200pt]{235.600pt}{0.400pt}}
\put(40,379){\makebox(0,0){Prob.}}
\put(650,21){\makebox(0,0){Total number of customers}}
\put(161.0,123.0){\rule[-0.200pt]{0.400pt}{123.341pt}}
\put(979,585){\makebox(0,0)[r]{Case (N, GD)}}
\put(999.0,585.0){\rule[-0.200pt]{24.090pt}{0.400pt}}
\put(161,507){\usebox{\plotpoint}}
\multiput(161.58,501.86)(0.497,-1.433){63}{\rule{0.120pt}{1.239pt}}
\multiput(160.17,504.43)(33.000,-91.428){2}{\rule{0.400pt}{0.620pt}}
\multiput(194.58,409.16)(0.497,-1.036){61}{\rule{0.120pt}{0.925pt}}
\multiput(193.17,411.08)(32.000,-64.080){2}{\rule{0.400pt}{0.463pt}}
\multiput(226.58,344.07)(0.497,-0.759){63}{\rule{0.120pt}{0.706pt}}
\multiput(225.17,345.53)(33.000,-48.535){2}{\rule{0.400pt}{0.353pt}}
\multiput(259.58,294.56)(0.497,-0.609){61}{\rule{0.120pt}{0.588pt}}
\multiput(258.17,295.78)(32.000,-37.781){2}{\rule{0.400pt}{0.294pt}}
\multiput(291.00,256.92)(0.549,-0.497){57}{\rule{0.540pt}{0.120pt}}
\multiput(291.00,257.17)(31.879,-30.000){2}{\rule{0.270pt}{0.400pt}}
\multiput(324.00,226.92)(0.719,-0.496){43}{\rule{0.674pt}{0.120pt}}
\multiput(324.00,227.17)(31.601,-23.000){2}{\rule{0.337pt}{0.400pt}}
\multiput(357.00,203.92)(0.895,-0.495){33}{\rule{0.811pt}{0.119pt}}
\multiput(357.00,204.17)(30.316,-18.000){2}{\rule{0.406pt}{0.400pt}}
\multiput(389.00,185.92)(1.113,-0.494){27}{\rule{0.980pt}{0.119pt}}
\multiput(389.00,186.17)(30.966,-15.000){2}{\rule{0.490pt}{0.400pt}}
\multiput(422.00,170.92)(1.486,-0.492){19}{\rule{1.264pt}{0.118pt}}
\multiput(422.00,171.17)(29.377,-11.000){2}{\rule{0.632pt}{0.400pt}}
\multiput(454.00,159.93)(2.145,-0.488){13}{\rule{1.750pt}{0.117pt}}
\multiput(454.00,160.17)(29.368,-8.000){2}{\rule{0.875pt}{0.400pt}}
\multiput(487.00,151.93)(2.476,-0.485){11}{\rule{1.986pt}{0.117pt}}
\multiput(487.00,152.17)(28.879,-7.000){2}{\rule{0.993pt}{0.400pt}}
\multiput(520.00,144.93)(3.493,-0.477){7}{\rule{2.660pt}{0.115pt}}
\multiput(520.00,145.17)(26.479,-5.000){2}{\rule{1.330pt}{0.400pt}}
\multiput(552.00,139.94)(4.722,-0.468){5}{\rule{3.400pt}{0.113pt}}
\multiput(552.00,140.17)(25.943,-4.000){2}{\rule{1.700pt}{0.400pt}}
\multiput(585.00,135.95)(6.937,-0.447){3}{\rule{4.367pt}{0.108pt}}
\multiput(585.00,136.17)(22.937,-3.000){2}{\rule{2.183pt}{0.400pt}}
\multiput(617.00,132.95)(7.160,-0.447){3}{\rule{4.500pt}{0.108pt}}
\multiput(617.00,133.17)(23.660,-3.000){2}{\rule{2.250pt}{0.400pt}}
\put(650,129.67){\rule{7.950pt}{0.400pt}}
\multiput(650.00,130.17)(16.500,-1.000){2}{\rule{3.975pt}{0.400pt}}
\put(683,128.17){\rule{6.500pt}{0.400pt}}
\multiput(683.00,129.17)(18.509,-2.000){2}{\rule{3.250pt}{0.400pt}}
\put(715,126.67){\rule{7.950pt}{0.400pt}}
\multiput(715.00,127.17)(16.500,-1.000){2}{\rule{3.975pt}{0.400pt}}
\put(748,125.67){\rule{7.709pt}{0.400pt}}
\multiput(748.00,126.17)(16.000,-1.000){2}{\rule{3.854pt}{0.400pt}}
\put(780,124.67){\rule{7.950pt}{0.400pt}}
\multiput(780.00,125.17)(16.500,-1.000){2}{\rule{3.975pt}{0.400pt}}
\put(846,123.67){\rule{7.709pt}{0.400pt}}
\multiput(846.00,124.17)(16.000,-1.000){2}{\rule{3.854pt}{0.400pt}}
\put(813.0,125.0){\rule[-0.200pt]{7.950pt}{0.400pt}}
\put(1009,122.67){\rule{7.709pt}{0.400pt}}
\multiput(1009.00,123.17)(16.000,-1.000){2}{\rule{3.854pt}{0.400pt}}
\put(878.0,124.0){\rule[-0.200pt]{31.558pt}{0.400pt}}
\put(1041.0,123.0){\rule[-0.200pt]{23.608pt}{0.400pt}}
\sbox{\plotpoint}{\rule[-0.400pt]{0.800pt}{0.800pt}}%
\put(979,524){\makebox(0,0)[r]{Case (N, GI)}}
\put(999.0,524.0){\rule[-0.400pt]{24.090pt}{0.800pt}}
\put(161,507){\usebox{\plotpoint}}
\multiput(162.41,494.30)(0.503,-1.815){59}{\rule{0.121pt}{3.061pt}}
\multiput(159.34,500.65)(33.000,-111.648){2}{\rule{0.800pt}{1.530pt}}
\multiput(195.41,379.04)(0.503,-1.394){57}{\rule{0.121pt}{2.400pt}}
\multiput(192.34,384.02)(32.000,-83.019){2}{\rule{0.800pt}{1.200pt}}
\multiput(227.41,294.03)(0.503,-0.932){59}{\rule{0.121pt}{1.679pt}}
\multiput(224.34,297.52)(33.000,-57.516){2}{\rule{0.800pt}{0.839pt}}
\multiput(260.41,235.02)(0.503,-0.625){57}{\rule{0.121pt}{1.200pt}}
\multiput(257.34,237.51)(32.000,-37.509){2}{\rule{0.800pt}{0.600pt}}
\multiput(291.00,198.09)(0.611,-0.504){47}{\rule{1.178pt}{0.121pt}}
\multiput(291.00,198.34)(30.555,-27.000){2}{\rule{0.589pt}{0.800pt}}
\multiput(324.00,171.09)(0.989,-0.507){27}{\rule{1.753pt}{0.122pt}}
\multiput(324.00,171.34)(29.362,-17.000){2}{\rule{0.876pt}{0.800pt}}
\multiput(357.00,154.08)(1.395,-0.511){17}{\rule{2.333pt}{0.123pt}}
\multiput(357.00,154.34)(27.157,-12.000){2}{\rule{1.167pt}{0.800pt}}
\multiput(389.00,142.08)(2.752,-0.526){7}{\rule{3.971pt}{0.127pt}}
\multiput(389.00,142.34)(24.757,-7.000){2}{\rule{1.986pt}{0.800pt}}
\multiput(422.00,135.06)(4.958,-0.560){3}{\rule{5.320pt}{0.135pt}}
\multiput(422.00,135.34)(20.958,-5.000){2}{\rule{2.660pt}{0.800pt}}
\put(454,128.84){\rule{7.950pt}{0.800pt}}
\multiput(454.00,130.34)(16.500,-3.000){2}{\rule{3.975pt}{0.800pt}}
\put(487,126.34){\rule{7.950pt}{0.800pt}}
\multiput(487.00,127.34)(16.500,-2.000){2}{\rule{3.975pt}{0.800pt}}
\put(520,124.84){\rule{7.709pt}{0.800pt}}
\multiput(520.00,125.34)(16.000,-1.000){2}{\rule{3.854pt}{0.800pt}}
\put(552,123.84){\rule{7.950pt}{0.800pt}}
\multiput(552.00,124.34)(16.500,-1.000){2}{\rule{3.975pt}{0.800pt}}
\put(585,122.84){\rule{7.709pt}{0.800pt}}
\multiput(585.00,123.34)(16.000,-1.000){2}{\rule{3.854pt}{0.800pt}}
\put(650,121.84){\rule{7.950pt}{0.800pt}}
\multiput(650.00,122.34)(16.500,-1.000){2}{\rule{3.975pt}{0.800pt}}
\put(617.0,124.0){\rule[-0.400pt]{7.950pt}{0.800pt}}
\put(683.0,123.0){\rule[-0.400pt]{109.850pt}{0.800pt}}
\end{picture}
\end{center}
\caption{Complementary distribution of total number of customers 
in Example 1.}
\label{figure-C}
\end{figure}

\begin{table}[tb]
\begin{center}
\caption{Joint queue length distribution 
$\protect\vc{p}(n_1, n_2)\protect\vc{e}$.}
(Upper rows for Case (P, GD) and lower rows for Case (P, GI) )\\
\label{table-14}
\begin{tabular}{@{}ccccc@{}}
\hline
\quad $n_1$
& \multicolumn{1}{c}{$0$}
& \multicolumn{1}{c}{$1$}
& \multicolumn{1}{c}{$2$}
& \multicolumn{1}{c}{$3$} 
\\ $n_2$
& \multicolumn{1}{c}{}
& \multicolumn{1}{c}{}
& \multicolumn{1}{c}{}
& \multicolumn{1}{c}{} 
\\
\hline
\down{2.5mm}{$0$} 
& \multicolumn{1}{c}{$2.500 \times 10^{-1}$}
& \multicolumn{1}{c}{$2.472 \times 10^{-2}$}
& \multicolumn{1}{c}{$8.593 \times 10^{-3}$}
& \multicolumn{1}{c}{$3.481 \times 10^{-3}$}
\\
& \multicolumn{1}{c}{$2.500 \times 10^{-1}$}
& \multicolumn{1}{c}{$4.501 \times 10^{-2}$}
& \multicolumn{1}{c}{$2.054 \times 10^{-2}$}
& \multicolumn{1}{c}{$9.224 \times 10^{-3}$}
\\
\cline{1-5}
\down{2.5mm}{$1$}
& \multicolumn{1}{c}{$6.530 \times 10^{-2}$}
& \multicolumn{1}{c}{$4.108 \times 10^{-2}$}
& \multicolumn{1}{c}{$2.193 \times 10^{-2}$}
& \multicolumn{1}{c}{$1.118 \times 10^{-2}$}
\\
& \multicolumn{1}{c}{$4.501 \times 10^{-2}$}
& \multicolumn{1}{c}{$4.108 \times 10^{-2}$}
& \multicolumn{1}{c}{$2.767 \times 10^{-2}$}
& \multicolumn{1}{c}{$1.629 \times 10^{-2}$}
\\
\cline{1-5}
\down{2.5mm}{$2$}
& \multicolumn{1}{c}{$3.249 \times 10^{-2}$}
& \multicolumn{1}{c}{$3.341 \times 10^{-2}$}
& \multicolumn{1}{c}{$2.444 \times 10^{-2}$}
& 
\\
& \multicolumn{1}{c}{$2.054 \times 10^{-2}$}
& \multicolumn{1}{c}{$2.767 \times 10^{-2}$}
& \multicolumn{1}{c}{$2.444 \times 10^{-2}$}
&  
\\
\cline{1-4}
\down{2.5mm}{$3$}
& \multicolumn{1}{c}{$1.497 \times 10^{-2}$}
& \multicolumn{1}{c}{$2.141 \times 10^{-2}$}
& 
& 
\\
& \multicolumn{1}{c}{$9.224 \times 10^{-3}$}
& \multicolumn{1}{c}{$1.629 \times 10^{-2}$}
& 
&  
\\
\cline{1-3}
\down{2.5mm}{$4$}
& \multicolumn{1}{c}{$6.630  \times 10^{-3}$}
& 
& 
& 
\\
& \multicolumn{1}{c}{$4.073 \times 10^{-3}$}
& 
&
&  
\\
\cline{1-2}
\end{tabular}
\end{center}
\end{table}

From Figures \ref{figure-B} and \ref{figure-C}, we observe that
class-dependent service times cause longer tails in the total queue
length distributions, in these specific examples. We shall explain
this phenomenon for Case N. In Case (N, GD), the conditional expected
amounts of work brought into the system per unit time given the state
of the underlying Markov chain are different, and they are given by
0.3 and 1.2, respectively. Thus in Case (N, GD), the system is
overloaded during a half of time. On the other hand, in Case (N, GI),
the conditional expected amount of work brought into the system per
unit time is fixed to be 0.75, regardless of the state of the
underlying Markov chain.  Therefore the distribution of the total
number of customers in Case (N, GD) has a longer tail than that in
Case (N, GI).

Next, we consider the expected total number $\E[N]$ of customers as a
function of the mean batch size $g$. Table \ref{table-11} shows
$\E[N]$ for the mean batch size $g=1$, 2, 3, 4, 5 and 10. We observe
that $\E[N]$ increases with the mean batch size $g$ in all cases. This
phenomenon comes from the fact that the deviation of the amount of
work brought into the system per unit time increases with $g$.  We
also observe that for a fixed $g$, the positive correlation in the two
streams leads to a larger $\E[N]$ in both Cases GD and GI, as
expected.
\begin{table}[tb]
\begin{center}
\caption{Expected total number of customers in Example 1.}
\label{table-11}
\begin{tabular}{@{\hspace*{3.6pt}}ccccccc@{\hspace*{3.6pt}}}
\hline
 Case & $g=1$ & $g=2$ & $g=3$ & $g=4$ & $g=5$ & $g=10$
\\
\hline
(P, GD) 
& 5.8760 & 9.9815 & 13.9356 & 17.8320 & 21.7001 & 40.8865
\\

\hline
(P, GI) 
& 5.8760 & 9.1466 & 12.2898 & 15.3793 & 18.4408 & 33.5873
\\
\hline
\hline
(I, GD) 
& 4.5417 & 8.5777 & 12.4865 & 16.3524 & 20.1987 & 39.3295
\\
\hline
(I, GI) 
& 4.0010 & 7.1857 & 10.2714 & 13.3219 & 16.3555 & 31.4326
\\
\hline
\hline
(N, GD) 
&3.2822 & 7.2033 & 11.0527 & 14.8822 & 18.7035 & 37.7739 
\\
\hline
(N, GI) 
& 2.2800 & 5.2800 & 
8.2800 & 11.2800 & 14.2800 & 29.2800
\\
\hline
\end{tabular}
\end{center}
\end{table}

\subsection{Number of customers in Example 2}

Table \ref{table-15} shows the expected total number $\E[N]$ of
customers for the mean batch size $g=1$, 2, 3, 4, 5. and 10. We first
examine the case of $g=1$. Contrary to Example 1, we observe that the
class-dependent service time (Case GD) decreases the expected total
number of customers in Cases I and N. This phenomenon can be explained
in a similar way to Example 1. For example, in Case (N, GD), the
conditional expected amount of work brought into the system per unit
time is fixed to be 0.8, regardless of the state of the underlying
Markov chain.  On the other hand, in Case (N, GI), the conditional
expected amounts of work brought into the system per unit time given 
the state of the underlying Markov chain are different, 
and they are given by 
1.28 and 0.32, respectively. Thus in Case (N, GI), the system is
overloaded during a half of time, so that $\E[N]$ in Case (N, GI) is
greater than that in Case (N, GD).

We observe that in any case, the expected total number of customer
increases with the mean batch size $g$, as in Example 1, and that
$\E[N]$ in Case GD eventually becomes greater than $\E[N]$ in Case GI.
We also observe that for a fixed $g$, the positive correlation in the
two streams leads to a larger $\E[N]$ in both Cases GD and GI, as in
Example 1.
\begin{table}[tb]
\begin{center}
\caption{Expected total number of customers in Example 2.}
\label{table-15}
\begin{tabular}{@{\hspace*{0pt}}c@{\hspace*{8.2pt}}cccccc@{\hspace*{0pt}}}
\hline
 Case & $g=1$ & $g=2$ & $g=3$ & $g=4$ & $g=5$ & $g=10$
\\
\hline
(P, GD) 
& 11.5019 & 17.7712 & 23.8347 & 29.8053 & 35.7261 & 65.0310
\\
\hline
(P, GI) 
& 11.5019 & 15.9366 & 20.1933 & 24.3657 & 28.4904 & 48.8117
\\
\hline
\hline
(I, GD) 
& 7.1517 & 13.3270 & 19.3007 & 25.2050 & 31.0760 & 60.2474
\\
\hline
(I, GI) 
& 8.7304 & 13.1052 & 17.3093 & 21.4407 & 25.5333 & 45.7640
\\
\hline
\hline
(N, GD) 
& 3.2168 & 9.0326 & 14.8425 & 20.6497 & 26.4551 & 55.4705
\\
\hline
(N, GI) 
& 6.0892 & 10.3399 & 14.4641 & 18.5407 & 22.5933 & 42.7206
\\
\hline
\end{tabular}
\end{center}
\end{table}

\section{Concluding Remarks}\label{sec:7}

We developed a numerically feasible procedure to compute the joint
queue length distribution in a FIFO single-server queue with multiple
batch Markovian arrival streams, under the assumption that service
time distributions of customers from respective arrival streams are
different and the batch size distributions follow discrete phase-type
distributions. We established several truncation and stopping criteria
to ensure the numerical accuracy in the final result.

Note, however, that the computation of the joint queue length
distribution is intensive by nature, especially when the number of
classes is large.  Even in such a case, the steady state distribution
of the total number of customers can be readily computed by modifying
our algorithm. For the sake of completeness, we show algorithm steps
for the total number of customers in Appendix. Note here that the
algorithm to compute $\vc{A}_k^{(\rm{T})}(n)$ in
(\ref{total_number_A_k(n)}) for the number of arrivals in a service
time can be used in the standard algorithm for the ordinary BMAP/GI/1
queue, too (see \cite{Luca91,Taki00}), because the BMAP/GI/1 queue is
considered as a special case of $K=1$ and the sequence of matrices for
the number of arrivals in a service time is essential for computing
the queue length distribution. To the best of our knowledge, however,
there is no work to consider the truncation and stopping criterion to
compute $\vc{A}_k^{(\rm{T})}(n)$ in the BMAP/GI/1 queue. Thus our
development also contributes to the standard algorithm for the
BMAP/GI/1 queue.

\appendix
\setcounter{section}{1}
\renewcommand{\thesection}{\Alph{section}}

\section*{Appendix: Algorithm for the Total Number of Customers}
\label{appendix-total}

We show a numerical algorithm to compute the steady state distribution
of the total number of customers, by modifying our algorithm for the
joint queue length distribution. We define $\vc{p}^{(\rm{T})}(n)$
($n=0,1,\ldots$) and $\vc{q}_k^{(\rm{T})}(n)$ ($k \in \K$,
$n=0,1,\ldots$) as
\[
\vc{p}^{(\rm{T})}(n) 
= \sum_{\scriptstyle \svc{n} \in \Z\atop \scriptstyle |\svc{n}| = n}
\vc{p}(\vc{n}),
\qquad
\vc{q}_k^{(\rm{T})}(n) 
= \sum_{\scriptstyle \svc{n} \in \Z \atop\scriptstyle|\svc{n}| = n} 
\vc{q}_k(\vc{n}),
\]%
\noindent
respectively. Corollary \ref{coro:a-04} is then reduced to
\begin{eqnarray*}
\vc{p}^{(\rm{T})}(0) 
&=& \sum_{k \in \K} \lambda_k \vc{q}_k^{(\rm{T})}(0)(-\vc{C})^{-1},
\\
\vc{p}^{(\rm{T})}(n)
&=& \sum_{k \in \K} 
\bigg[ \lambda_k 
\left( \vc{q}_k^{(\rm{T})}(n) - \vc{q}_k^{(\rm{T})}(n-1) \right)
\nonumber
\\
&& \qquad {} + \sum_{m=1}^{n} 
\vc{p}^{(\rm{T})}(n-m) \vc{D}_k(m) 
\bigg] (-\vc{C})^{-1}, 
\quad n =1,2,\dots.
\end{eqnarray*}%
\noindent
Further, under Assumption \ref{ass:a-01}, Theorem \ref{thm:a-01} is
reduced to
\begin{eqnarray*}
\vc{q}_k^{(\rm{T})}(n) 
&=& {1 \over \lambda_k}
\sum_{\scriptstyle m_1+m_2+m_3
\atop \scriptstyle + m_4 = n }
\vc{v}_k^{(\rm{T})}(m_1) 
[ \vc{\alpha}_k \otimes \vc{A}_k^{(\rm{T})}(m_2)]
\\ 
&& \qquad {} \cdot 
\vc{\Gamma}_k^{(\rm{T})}(m_3)
\left[ \{ \vc{P}_k^{m_4} (\vc{I} - \vc{P}_k) \vc{e} \} 
\otimes \vc{I}(M) \right],
\end{eqnarray*}%
\noindent
where 
\begin{eqnarray}
\vc{A}_k^{(\rm{T})}(n) 
&=& \sum_{\scriptstyle \svc{n} \in \Z\atop \scriptstyle |\svc{n}| = n}
\vc{A}_k(\vc{n}),
\label{total_number_A_k(n)}
\\
\vc{v}_k^{(\rm{T})}(n) 
&=& \sum_{\scriptstyle \svc{n} \in \Z\atop \scriptstyle |\svc{n}| = n} 
\vc{v}_k(\vc{n}),
\qquad
\vc{\Gamma}_k^{(\rm{T})}(n) 
= \sum_{\scriptstyle \svc{n} \in \Z\atop \scriptstyle |\svc{n}| = n}
\vc{\Gamma}_k(\vc{n}).
\nonumber
\end{eqnarray}%
\noindent
Thus the $\vc{p}^{(\rm{T})}(n)$ is obtained if we compute the
$\vc{A}_k^{(\rm{T})}(n)$, the $\vc{v}_k^{(\rm{T})}(n)$ and the
$\vc{\Gamma}_k^{(\rm{T})}(n)$.

Note here that $\vc{A}_k^{(\rm{T})}(n)$, $\vc{v}_k^{(\rm{T})}(n)$ and
$\vc{\Gamma}_k^{(\rm{T})}(n)$ satisfy
\begin{eqnarray*}
\sum_{n=0}^{\infty} z^n
\vc{A}_k^{(\rm{T})}(n) 
&=& \int_0^{\infty} dH_k(x)
\exp\left[ \left(\vc{C} + \sum_{k \in \K} \vc{D}_k^{\ast}(z) 
\right) x\right],
\\
\sum_{n=0}^{\infty} z^n 
\vc{v}_k^{(\rm{T})}(n) 
&=& \int_0^{\infty} d\vc{v}(x) \vc{D}_k 
\exp\left[ \left(\vc{C} + \sum_{k \in \K} \vc{D}_k^{\ast}(z) 
\right) x\right],
\\
\sum_{n=0}^{\infty} z^n \vc{\Gamma}_k^{(\rm{T})}(n) 
&=& \left[ \vc{I} - \vc{P}_k \otimes \int_0^{\infty} dH_k(x)  
\exp\left[ \left(\vc{C} + \sum_{k \in \K} \vc{D}_k^{\ast}(z) 
\right) x\right]
\right]^{-1},
\end{eqnarray*}%
\noindent
respectively. Thus $\vc{A}_k^{(\rm{T})}(n)$ and
$\vc{v}_k^{(\rm{T})}(n)$ can be written to be
\begin{eqnarray*}
\vc{A}_k^{(\rm{T})}(n) 
&=& \sum_{m=0}^{\infty}\gamma_k^{(m)}(\theta)\vc{F}_m^{(\rm{T})}(n),
\\
\vc{v}_k^{(\rm{T})}(n) 
&=& \sum_{m=0}^{\infty}\vc{v}^{(m)}(\theta)\vc{D}_k\vc{F}_m^{(\rm{T})}(n),
\end{eqnarray*}%
\noindent
respectively, where $\vc{F}_m^{(\rm{T})}(n)$ denotes an $M \times M$
matrix which satisfies
\[
\sum_{n=0}^{\infty} z^n \vc{F}_m^{(\rm{T})}(n) 
= \left[\vc{I} + \theta^{-1}
\left(\vc{C} + \sum_{k \in \K} \vc{D}_k^{\ast}(z) \right) \right]^m.
\]%
\noindent
Further the $\vc{\Gamma}_k^{(\rm{T})}(n)$ ($k \in \K$, $n \ge 0$) is
determined by the following recursion:
\[
\vc{\Gamma}_k^{(\rm{T})}(0) 
= \left[\vc{I} - \vc{P}_k \otimes 
\vc{A}_k^{(\rm{T})}(0) \right]^{-1}, 
\]%
\noindent
and for $n =1,2,\dots$,
\[
\vc{\Gamma}_k^{(\rm{T})}(n) 
= \sum_{l=1}^n \vc{\Gamma}_k^{(\rm{T})}(n-l)
\left[\vc{P}_k \otimes \vc{A}_k^{(\rm{T})}(l) 
\right]\vc{\Gamma}_k^{(\rm{T})}(0).
\]%
\noindent
Thus we can compute $\vc{A}_k^{(\rm{T})}(n)$, $\vc{v}_k^{(\rm{T})}(n)$
and $\vc{\Gamma}_k^{(\rm{T})}(n)$ by replacing Steps 3 and 4 with the
followings.
\medskip

\noindent
{\bf Step 3. } Compute $\breve{\vc{A}}_k^{(\rm{T})}(n)$ and
$\breve{\vc{v}}_k^{(\rm{T})}(n)$ by the following procedure, where the
initial values of $\breve{\vc{A}}_k^{(\rm{T})}(n)$ and
$\breve{\vc{v}}_k^{(\rm{T})}(n)$ ($n \ge 0$) are assumed to be
$\vc{O}$ and $\vc{0}$, respectively.

\begin{narrow}
{\bf Step (3--a). } Set $\breve{\vc{F}}_0^{(\rm{T})}(0) = \vc{I}$ and
$n_F^{(0)} = 0$.  Also set
\[
\breve{\vc{A}}_k^{(\rm{T})}(0) = \gamma_k^{(0)}(\theta)\vc{I},
\quad
\breve{\vc{v}}_k^{(\rm{T})}(0) = \vc{v}^{(0)}(\theta)\vc{D}_k,
\quad \forall k \in \K.
\]%
\end{narrow}

\begin{narrow}
{\bf Step (3--b). } Set $n_F^{(1)} = \max_{k \in \K} n_g(k)$ and
$m=1$, and compute $\breve{\vc{F}}_1^{(\rm{T})}(n)$ by the following
recursion:
\[
\breve{\vc{F}}_1^{(\rm{T})}(0) 
= \vc{I} + \theta^{-1}\vc{C},
\]%
\noindent
and for $n=1,2,\dots,n_F^{(1)}$,
\[
\breve{\vc{F}}_1^{(\rm{T})}(n)
=
\theta^{-1} \sum_{k \in \K} U(n_g(k) - n)g_k(n) \vc{D}_k.
\]%
\end{narrow}

\begin{narrow}
{\bf Step (3--c). } For each $k \in \K$, if $m \leq m_{\gamma}(k)$,
add $\gamma_k^{(m)}(\theta) \breve{\vc{F}}_m^{(\rm{T})}(n)$ to
$\breve{\vc{A}}_k^{(\rm{T})}(n)$ for all $n \leq n_F^{(m)}$.  Also,
for each $k \in \K$, if $m \leq m_{v}(k)$, add
$\vc{v}^{(m)}(\theta)\vc{D}_k \breve{\vc{F}}_m^{(\rm{T})}(n)$ to
$\breve{\vc{v}}_k^{(\rm{T})}(n)$ for all $n \leq n_F^{(m)}$.
\end{narrow}

\begin{narrow}
{\bf Step (3--d). }
If $m \geq m_{\max}$, stop computing, and otherwise, 
add one to $m$ and go to Step (3--e).
\end{narrow}

\begin{narrow}
{\bf Step (3--e). } For each $n=0,1,\dots$, compute
$\breve{\vc{F}}_m^{(\rm{T})}(n)$ by
\begin{eqnarray*}
\breve{\vc{F}}_m^{(\rm{T})}(n)
&=& U\left(n_F^{(m-1)} - n\right)
\breve{\vc{F}}_{m-1}^{(\rm{T})}(n)(\vc{I} + \theta^{-1}\vc{C}) 
\nonumber
\\
& & \quad {} + 
\sum_{l=1}^{\min(n,\ n_F^{(1)})}
U\left(n_F^{(m-1)} - n  + l  \right)
\breve{\vc{F}}_{m-1}^{(\rm{T})}(n - l) 
\breve{\vc{F}}_1^{(\rm{T})}(l),
\end{eqnarray*}%
\noindent
until $\breve{\vc{F}}_m^{(\rm{T})}(n)$'s satisfy $\sum_{n \le
n^{\ast}} \breve{\vc{F}}_m^{(\rm{T})}(n) \vc{e} > \left(1 -
\varepsilon_F \right)^m \vc{e}$ for some $n^{\ast}$.  Let $n_F^{(m)} =
n^{\ast}$ and go to Step (3--c).
\end{narrow}

\noindent
{\bf Step 4. }
Set 
\[
n_A(k)
= \max \left( n_F^{(m)} ; m=0,1,\dots,m_{\gamma}(k) \right),
\]%
\noindent
and for each $k \in \K$, compute 
$\breve{\vc{\Gamma}}_k^{(\rm{T})}(n)$ by the following recursion:
\[
\breve{\vc{\Gamma}}_k^{(\rm{T})}(0) 
= \left[\vc{I} - \vc{P}_k \otimes 
\breve{\vc{A}}_k^{(\rm{T})}(0) \right]^{-1}, 
\]%
\noindent
and for $n=1,2,\dots$,
\begin{eqnarray*}
\breve{\vc{\Gamma}}_k^{(\rm{T})}(n) 
&=& \sum_{l=1}^n 
U\left(n_A(k) - l\right) 
\breve{\vc{\Gamma}}_k^{(\rm{T})}(n-l)
\left[\vc{P}_k \otimes \breve{\vc{A}}_k^{(\rm{T})}(l) 
\right]\breve{\vc{\Gamma}}_k^{(\rm{T})}(0),
\end{eqnarray*}%
\noindent
until $\breve{\vc{\Gamma}}_k^{(\rm{T})}(n)$'s satisfy
\begin{eqnarray*}
\sum_{n=0}^{n_{\Gamma}(k)} 
\breve{\vc{\Gamma}}_k^{(\rm{T})}(n)\vc{e}
> \left\{ \left(\vc{I} - \vc{P}_k\right)^{-1} \vc{e}(M_k) \right\}
\otimes \vc{e}(M)
-
\varepsilon 
\left\{ \left(\vc{I} - \vc{P}_k\right)^{-2} \vc{P}_k \vc{e}(M_k) \right\} 
\otimes \vc{e}(M),
\end{eqnarray*}%
\noindent
for some integer $n_{\Gamma}(k)$.

\begin{rem}
The above algorithm ensures that 
\[
\sum_{n=0}^{n_A(k)}
\breve{\vc{A}}_k^{(\rm{T})}(n) \vc{e}
> (1 - \varepsilon)\vc{e},
\quad 
\sum_{n=0}^{n_v(k)} 
\breve{\vc{v}}_k^{(\rm{T})}(n) \vc{e}
> (1 - \varepsilon)\lambda_k^{(\rm{B})},
\]%
\noindent
respectively, where $n_v(k)$ is given by
\[
n_v(k) 
= \max \left( n_F^{(m)} ; m=0,1,\dots,m_v(k) \right).
\]%
\noindent
Further $\breve{\vc{\Gamma}}_k^{(\rm{T})}(n)$ satisfies
\[
(\vc{\alpha}_k \otimes \vc{\pi} )
\sum_{n=0}^{n_{\Gamma}(k)} \breve{\vc{\Gamma}}_k^{(\rm{T})}(n)\vc{e}
> \E[G_k] - {1 \over 2} \E[G_k(G_k-1)] \varepsilon.
\]%
\end{rem}

\begin{small}

\end{small}
\end{document}